\documentclass[a4paper,11pt]{article}
\usepackage[english]{babel}
\usepackage[utf8]{inputenc}
\usepackage[T1]{fontenc}
\usepackage{textcomp} 
\usepackage{amsmath}
\usepackage{enumerate}
\usepackage{amsfonts}
\usepackage{amsthm}
\usepackage{mathrsfs}
\usepackage{amssymb}
\usepackage{csquotes}
\usepackage{graphicx}
\usepackage{bbm}
\usepackage{mathabx}
\usepackage[mathscr]{euscript}
\DeclareSymbolFont{rsfs}{U}{rsfs}{m}{n}
\DeclareSymbolFontAlphabet{\mathscrsfs}{rsfs}

\newtheorem{prop}{Proposition}
\newtheorem{lemma}{Lemma}
\newtheorem{definition}{Definition}
\newtheorem{corollary}{Corollary}
\newtheorem{theorem}{Theorem}

\def\real{{\mathord{{\rm I\kern-2.8pt R}}}}        
\def\inte{{\mathord{{\rm I\kern-2.8pt N}}}}

\def\sZZ{{\rm Z\kern-2.8ptem{}Z}}

\def\z{{\mathchoice
		{\sZZ}
		{\sZZ}
		{\rm Z\kern-0.30em{}Z}
		{\rm Z\kern-0.25em{}Z} }}
\def\sQQ{{\kern 0.27em \vrule height1.45ex width0.03em depth0em
		\kern-0.30em \rm Q}}
\def\qu{{\mathchoice
		{\sQQ}
		{\sQQ}
		{\kern 0.225em \vrule height1.05ex width0.025em depth0em \kern-0.25em \rm Q}
		{\kern 0.180em \vrule height0.78ex width0.020em depth0em \kern-0.20em \rm Q}
}}
\def\sCC{{\kern 0.27em \vrule height1.45ex width0.03em depth0em
		\kern-0.30em \rm C}}
\def\complex{{\mathchoice
		{\sCC}
		{\sCC}
		{\kern 0.225em \vrule height1.05ex width0.025em depth0em \kern-0.25em \rm C}
		{\kern 0.180em \vrule height0.78ex width0.020em depth0em \kern-0.20em \rm C}
}}

\newcommand{\R}{\mathbb{R}}

\newcommand{\ba}{\begin{array}}
	\newcommand{\ea}{\end{array}}
\newcommand{\be}{\begin{equation}}
	\newcommand{\ee}{\end{equation}}
\newcommand{\bea}{\begin{eqnarray}}
	\newcommand{\eea}{\end{eqnarray}}
\newcommand{\beaa}{\begin{eqnarray*}}
	\newcommand{\eeaa}{\end{eqnarray*}}

\def\z{\zeta}

\font\tenmath=msbm10 \font\sevenmath=msbm7 \font\fivemath=msbm5
\newfam\mathfam \textfont\mathfam=\tenmath
\scriptfont\mathfam=\sevenmath \scriptscriptfont\mathfam=\fivemath
\def\math{\fam\mathfam}

\def \={{\buildrel {\rm (law)} \over =}}

\def \R{{\math R}}
\def \Z{{\math Z}}

\def\qed{ \hfill \vrule width.25cm height.25cm depth0cm\smallskip}

\newcommand{\basa}{\begin{assumption}}
	\newcommand{\easa}{\end{assumption}}

\newcommand{\bas}{\begin{assum}}
	\newcommand{\eas}{\end{assum}}

\newcommand{\ignore}[1]{}

\begin{document}
	
	\renewcommand{\thefootnote}{\fnsymbol{footnote}}
	
	\renewcommand{\thefootnote}{\fnsymbol{footnote}}

	\title{Modified wavelet variation for the Hermite processes}
	
	\author{L. Loosveldt\footnote{ Universit\'e de Li\`ege, D\'epartement de Math\'ematique -- zone Polytech 1, 12 all\'ee de la D\'ecouverte, B\^at. B37, B-4000 Li\`ege. l.loosveldt@uliege.be}, C.A. Tudor\footnote{Université de Lille, 
CNRS, UMR 8524 - Laboratoire Paul Painlev\'e, F-59000 Lille. ciprian.tudor@univ-lille.fr} }
	
	\maketitle
	
\begin{abstract}
We define an asymptotically normal wavelet-based strongly consistent estimator for the Hurst parameter of any Hermite processes. This estimator is obtained by considering a modified wavelet variation in which coefficients are wisely chosen to be, up to negligeable remainders, independent. We use Stein-Malliavin calculus to prove that this wavelet variation satisfies a multidimensional Central Limit Theorem, with an explicit bound for the Wasserstein distance.
\end{abstract}

	
\section{Introduction}

A stochastic process $\{X_t\}_{t \geq 0}$ is said to be \textit{$H$-self-similar} if, for all $a >0$, the processes $\{X_{at}\}_{t \geq 0}$ and $\{a^H X_t\}_{t \geq 0}$ have the same finite dimensional distributions. This property is predominant in real life applications, for instance, in astronomy \cite{Lomax2018,Pan2007}, biology \cite{Audit2002,Ecksteina}, climatology \cite{npg-14-723-2007,Franzke2020}, hydrology \cite{McLeod1978PreservationOT,Taqqu1978}, image processing \cite{MWEMA202013}, internet traffic modelling \cite{VivekKumarSehgal_2019,Willinger1995}, mathematical finance \cite{Fauth2016,Mandelbrot2012,Odonkor2019,STOYANOV2019} and physics \cite{Metzler2014,Perrin1910}. We also refer to the three monographs \cite{Beran1994,Doukhan2003,Embrechts2002} for an insight on numerous results, applications and methodologies concerning self-similar processes.

Without any doubt, fractional Brownian motion is the most popular self-similar stochastic process. It was first introduced by Kolmogorov, in the paper \cite{MR0003441} from 1940, to define ``Gaussian spirals'' in Hilbert spaces. The name ``fractional Brownian motion'' is used since the article \cite{MR242239} from Mandelbrot and Van Ness, where the first systematic study of this process was carried out. Indeed, it is itself a generalization of the Brownian motion, defined by the botanist Robert Brown to describe the movements of pollen grains of the plant Clarkia Pulchella suspended in the water \cite{brown}. For any $H \in (0,1)$, the fractional Brownian motion of Hurst parameter $H$ on a probability space $(\Omega,\mathcal{F},\mathbf{P})$ is the unique centred Gaussian process $\{B^H_t\}_{t \geq 0}$ with covariance function
\[ \mathbf{E}[B^H_tB^H_S] = \frac{1}{2}(|t|^{2H}+|s|^{2H}-|t-s|^{2H}), \mbox{ for every } s, t\geq 0. \]
The case $H=1/2$ corresponds to the usual Brownian motion. In other words, the fractional Brownian motion $\{B^H_t\}_{t \geq 0}$  is the only Gaussian process with stationary increments which is $H$-self-similar. This property justifies the fact that the fractional Brownian is greatly used to model natural phenomena \cite{Ecksteina,Franzke2020,Lomax2018,Metzler2014,Pan2007}. 

In some applications, it happens that the self-similarity and the stationnarity of increments are desirable properties for a stochastic process, while Gaussianity is not a reasonable assumption \cite{Fauth2016,STOYANOV2019,Taqqu1978,VivekKumarSehgal_2019,MR2212691}. In these contexts, good candidates for simulation are given by Hermite processes. In the sequel, for all integer number $q\geq 1$, $ I_{q}$ denotes the multiple stochastic integral of order $q$ with respect to the two-sided Brownian motion $\{B_y\}_{y \in \mathbb{R}}$, see the Appendix for a precise definition. Let $H \in \left( \frac{1}{2}, 1\right)$, the Hermite process of order $q$ and Hurst parameter $H$ is the stochastic process $\{ Z^{(q, H)} _{t}\}_{t\geq 0}$ defined, for every $t \geq 0$, by
\begin{equation}\label{h1}
	Z^{(q, H)}_{t}= I_{q} (L^{(q, H)}_{t}),
\end{equation}
where the kernel function $ L^{(q, H)}_{t}$ is given, for all $ y_{1},..., y_{q} \in \mathbb{R}$, by
\begin{equation}
	\label{M}
	L^{(q, H)}_{t} (y_{1},..., y_{q})=c_{q,H} \int_{\mathbb{R}} (u-y_{1})_{+} ^{-\left( \frac{1}{2}+\frac{1-H}{q}\right)}\ldots (u-y_{q})_{+} ^{-\left( \frac{1}{2}+\frac{1-H}{q}\right)}du,
\end{equation}
with $c_{q,H}$ a strictly positive constant which is chosen such that, for all $t \geq 0$, $ \mathbf{E} [(Z^{(q, H)}_{t}) ^{2}] =q! \Vert L^{(q, H)}_{t} \Vert ^{2} _{ L^2(\R ^{q})}=t ^{2H}$. It is a $H$-self-similar stochastic process with stationary increments. Its exhibits long-range dependence and its sample paths are H\"older continuous of order $\delta$ for every $\delta \in (0, H)$. Hermite processes first appeared as limit of partial sums of correlated random variables, in the so-called Non-Central Limit Theorem \cite{Dobrushin1979,Taqqu1974/75,Taqqu1979}. The class of Hermite processes contains the fractional Brownian motion which is obtained for $q=1$ and it is the only Gaussian process in this class. It also contained the Rosenblatt process (obtained for $q=2$). We refer to the recent monograph \cite{tudor2023} for a concise presentation of these processes and their stochastic analysis. The results present on this paper do not depend on the exact value of $q$ and $H$. For this reason, to ease the notation, we omit the indices and we write $\{Z_t\}_{t \geq 0}$ for a Hermite process as well as $L_t$ for its associated kernel.

Generally speaking, while working with self-similar processes, a question of great interest is to estimate the parameter $H$ here over. It is essential because this parameter, governs the main properties of the analysed process, such as the long-range dependence and the Hölder regularity. Therefore, in the Hermite case, many authors have proposed statistical estimators from various perspectives, such as variations \cite{Chronopoulou2011,Coeurjolly2001,Tudor2009}, wavelet analysis \cite{BaTu,Bardet2002,Clausel2014,Flandrin1992} or least-squares methods \cite{Nourdin2019}.

Unfortunately, all the estimators presented in the publications listed here over failed to be asymptotically normal as soon as $q \geq 2$. Also, in the case $q=1$, the asymptotic normality only holds when $H \in (0,\frac{3}{4}]$. Note that, in this last Gaussian case, the authors in \cite{Istas1997} propose an estimator based on higher order increments of the fractional Brownian motion and obtain the asymptotic normality of this estimator. But, as remarked in \cite{Chronopoulou2009}, this strategy is not applicable for $q \geq 2$. Of course, this non-Gaussian behavior is a serious drawback for statistical estimation of the Hurst parameter $H$, as the limit distribution is a Rosenblatt one, much less desirable for practical purposes.

But, very recently, in \cite{AT}, Ayache and Tudor define a new estimator for $H$, based on a modified quadratic variation, which is asymptotically normal, even in the case $q \geq 2$. The construction of this estimator is based on an idea from \cite{A} where increments of Hermite processes are split in two parts: one which satisfies some independence properties while the second part is clearly dominated, in $L^2(\Omega)$-norm, by the first one. These observations help to select ``good'' increments to compute quadratic variations which are then collected in a simple modification of the estimators previously introduced in \cite{Chronopoulou2011,Tudor2009}.

Nevertheless, in practice, estimators for the Hurst parameter based on wavelet analysis provide numerous advantages, compared to the ones obtained with other strategies. Indeed, wavelet-based estimators are defined using a log-regression of wavelet coefficients at various scales, see Section \ref{sec:wavelet} for the definition, and a great literature is available concerning goodness-of-fit tests for such models \cite{Bardet2010,Bardet2008,Bardet2002}. Moreover, such estimators are numerically efficient thanks to the Mallat's algorithm for computing wavelet coefficients \cite{mallat89,Mallat1998}. Also, if the wavelet have enough vanishing moment, see equation \eqref{mom1} below, wavelet-based estimators are not sensible to polynomials trends and are thus very robust.

In the paper \cite{Daw2022}, inspired by \cite{A}, Daw and Loosveldt adapted the ``splitting method'' in the context of wavelet coefficients for the Rosenblatt process (aka the Hermite process of order $2$). In this context, they manage to express any wavelet coefficient as the sum of two random variables, one which satisfies some independence property and one which is negligible.

The aim of this paper is to bring together the ideas from \cite{AT} and \cite{Daw2022} and to propose a new strongly consistent estimator for the Hurst parameter of any Hermite processes. It belongs to the ``family'' of the estimators introduced in \cite{BaTu,Bardet2002,Clausel2014,Flandrin1992}, in the sense that it relies on a wavelet variation. But the wavelet coefficients used in this ``modified'' variation are precisely chosen such that this new estimator is asymptotically normal. In short, this estimator is particularly well-suited for applications as it is an asymptotically normal wavelet-based strongly consistent estimator. Let us note that, up to our knowledge, it is anyway the first time that a wavelet-based estimator is defined for Hermite processes of any orders. Indeed, the papers \cite{BaTu,Bardet2002,Clausel2014,Flandrin1992} only deal with the fractional Brownian motion and the Rosenblatt process, thus the Hermite processes of order one and two.

In Section \ref{sec:wavelet}, we recall some basic facts about wavelet analysis on which relies the rest of our paper. Then we extend to any Hermite process the results from \cite{Daw2022} concerning the decomposition of wavelet coefficients.

In Section \ref{sec:modifiedvar}, we introduce the modified wavelet variation which is the key object to define later our estimator. First, we present the special wavelet coefficients that are used in the variation. Then, we show that a part of this variation is clearly negligible compared to the other part. It is a consequence of the decomposition of wavelet coefficients provided in Section \ref{sec:wavelet}. Finally, we prove a multidimensional Central Limit Theorem (CLT) for the modified wavelet variation. Note that the use of a multidimensional CLT is a major difference with \cite{AT}, where only a ``classic'' CLT is used. It is a consequence of the fact that our estimator relies on wavelet coefficients and not directly on increments of the process. Indeed, as the reader can see all along this paper, the use of a wavelet $\Psi$ make appeared a constant $C_\Psi(H)$ in all our computations, see equation \eqref{cpsi} below. As this constant depends, in a non-trivial way, on the value of the Hurst parameter $H$, it must be eliminated in the definition of the estimator. Then, the standard strategy from \cite{BaTu,Bardet2002} consists in performing a $\log$-regression of the wavelet variation onto several scales. Thus, we need to insure the asymptotic normality for a vector of wavelet variation instead of an unique variation based on increments.

In order to obtain an estimator which can be numerically computed, in Section \ref{sec:discret} we study a Riemann approximation of the wavelet variation defined in Section \ref{sec:modifiedvar}. The main result there states that this approximation converges, in $L^1(\Omega)$, to the wavelet variation and thus also satisfies a multidimensional CLT.

Finally, in Section \ref{sec:estima}, we use a $\log$-regression to define an estimator for the Hurst parameter $H$. We show that this estimator is strongly consistent and asymptotically normal.

Main facts of stochastic analysis, in particular Malliavin-Stein method, are recalled, for the reader convenience, in the Appendix.

\section{The wavelet coefficients and their decomposition}\label{sec:wavelet}

Generally speaking, a wavelet is a smooth function $\Psi$ such that the set 
\[ \{ 2^\frac{j}{2}\Psi(2^j \cdot-k) \, : \, (j,k) \in \Z^2 \}\]
is an orthonormal basis of $L^2(\R)$, see \cite{Daubechies1992,Mallat1998,Meyer1992} for a very complete view on this subject. Morally speaking, it is not surprising that such bases, constructed by scaling modifications of an unique mother wavelet $\Psi$ are particularly well-adapted to study self-similar stochastic processes, as this definition somehow means that the process is invariant by scaling changes. In particular, in the context of Hermite processes and their generalizations, wavelets are used, for instance, for expansion and approximation \cite{Ayache2020,ayachehamonierloosveldt,Ayache2005,Meyer1999,Pipiras2004}, estimation \cite{BaTu,Bardet2002} and precise study of the pointwise regularity \cite{Daw2022,Esser2022,esserloosveldt2}.

All along this paper, we work with a continuously differentiable mother wavelet $\Psi: \mathbb{R} \to \mathbb{R}$ with support included in the unit interval $[0,1]$ such as the one constructed in \cite {Daubechies1988,Daubechies1992}. We say that such a wavelet $\Psi$ has $Q \geq 1$ vanishing moments if, for all $0 \leq p < Q$, 
\begin{equation}
	\label{mom1}
	\int_{\mathbb{R}} x^p \Psi (x) dx =0.
\end{equation}
In fact, in this paper, we only need that $\Psi$ has one vanishing moment. But, in practice, authors generally work with higher $Q$, for numerical purposes. Let us fix $a>0$. We define, for any integer number $k\geq 0$, the wavelet coefficient associated to the Hermite process $\{Z_t\}_{t \geq 0}$ by 
\begin{equation}
	\label{coef}
	c(a,k)= \sqrt{a} \int_{\mathbb{R}} \Psi (x) Z _{ a (x+k)} dx.
\end{equation}
Due to the scaling property,  the stationarity of the increments of the Hermite process and the vanishing moment of $\Psi$, we notice
\begin{eqnarray}\label{2d-1}
	c(a,k)= \sqrt{q} \int_{\mathbb{R}} \Psi (x) ( Z_{a(x+k)}- Z_{ak})dx = ^{(d)} \sqrt{a} \int_{\mathbb{R}}\Psi (x) Z_{ax}dx =^{(d)} a ^{H+\frac{1}{2}} c(1,0),
\end{eqnarray}
where we denoted by ``$= ^{(d)}$'' the equality in distribution. In particular, we have
\begin{equation}\label{30n-3}
	\mathbf{E}[ c(a,k) ^{2}]= a ^{2H+1} 	\mathbf{E}[c(1,0) ^{2}] = a ^{2H+1} C_{\Psi}(H),
\end{equation}
with 
\begin{equation}\label{cpsi}
	C_{\Psi }(H) =-\frac{1}{2}\int_{\mathbb{R}} \int_{\mathbb{R}} dxdy \Psi (x) \Psi (y) \vert x-y\vert ^{2H}.
\end{equation}

By using (\ref{h1}) and the assumption (\ref{mom1}),  we can express the wavelet coefficient as a multiple stochastic integral. Indeed, we write, via a Fubini argument,  

\begin{eqnarray*}
	c(a,k) &=&\sqrt{a} \int_{\mathbb{R}}\Psi (x) ( Z_{a(x+k)}-Z_{ak})dx = \sqrt{a} \int_{\mathbb{R}}\Psi (x) I_{q} \left( L_{a(x+k)}-L_{ak}\right) dx \\
	&=&I_{q} \left( \sqrt{a} \int_{\mathbb{R}}\Psi (x)  \left( L_{a(x+k)}-L_{ak}\right) dx\right)\\
	&=& I_{q} \left( \sqrt{a} \int_{\mathbb{R}} dz \Psi (x) \int_{ak} ^{a(x+k)} f_{u} du \right),
\end{eqnarray*}
where $L$ is the kernel of the Hermite process given by (\ref{M}) and
\begin{equation*}
	f_{u} (y_{1},..., y_{q}) = c_{q,H}(u-y_{1})_{+} ^{-\left( \frac{1}{2}+\frac{1-H}{q}\right)}\ldots (u-y_{q})_{+} ^{-\left( \frac{1}{2}+\frac{1-H}{q}\right)},
\end{equation*}
for every $y_{1},..., y_{q} \in \mathbb{R}$, with $c_{q,H}$ the constant in \eqref{M}.  Since  $\int_{ak} ^{a(x+k)} f_{u}(y_{1},..., y_{q}) du$  vanishes if there exists $i=1,..,q $ such that $ y_{i} \geq a(k+1)$ (recall that the support of $\Psi $ is included in $[0,1]$), we can also write 

\begin{equation}
	\label{coef1}
	c(a,k) = I_{q} \left(  \sqrt{a} \int_{\mathbb{R}}dx \Psi (x) \mathbbm{1}_{ (-\infty, a(k+1))}^{\otimes q} \int_{ak} ^{a(x+k)} f_{u} du \right).
\end{equation}

Let $M>0$. We decompose the coefficient $c(a,k) $ as follows 
\begin{equation*}
	c(a,k)= \widetilde{c(a,k, M)}+\widecheck{c(a,k, M)}
\end{equation*}
where
\begin{equation}\label{c1}
	\widetilde{c(a,k, M)}= I_{q} \left(  \sqrt{a} \int_{\mathbb{R}}dx \Psi (x) \mathbbm{1}_{ (a(k-M), a(k+1))}^{\otimes q} \int_{ak} ^{a(x+k)} f_{u} du \right)
\end{equation}
and
\begin{equation}\label{c2}
	\widecheck{c(a,k, M)}= I_{q} \left(  \sqrt{a} \int_{\mathbb{R}}dx \Psi (x) \left( \mathbbm{1}_{ (-\infty, a(k+1))^{q}\setminus  (a(k-M), a(k+1))^{ q}}\right)\int_{ak} ^{a(x+k)} f_{u} du \right).
\end{equation}
This decomposition will play a crucial role for the development of our arguments. The key fact is that, for suitable choices of the parameters $a,M>0$, the random variables $\widetilde{c(a,k, M)}$ and  $\widetilde{c(a,j, M)}$ will become independent when $k\not=j$ (this happens once the intervals $ (a(k-M), a(k+1))$ and $(a(j-M), a(j+1))$ become disjoints). On the other hand, the terms $\widecheck{c(a,k,M)}$ will be negligible, for the  $ L^{2}(\Omega)$-norm, with respect to their tilde counterparts. 

Let us first evaluate the squared mean of $\widecheck{c(a,k,M)}$. We denoted by $C(q,H)$ a strictly positive constant depending only on $q,H$.  A similar notation will be used for constants depending on other parameters.

\begin{lemma}\label{ll1}
	Let $a,M>0$ and let $k$ be a positive integer. We have 
	\begin{equation*}
		\mathbf{E} \Big[ \widecheck{c (a, k,M)} ^{2} \Big] \leq  C(q,H) a ^{2H+1} (M+1)^{\frac{2H-2}{q}}.
	\end{equation*}
\end{lemma}
\noindent {\bf Proof: }We first notice that, for every $a, k, M$, we have, by the same arguments as in the proofs of Lemma \ref{ll2} below,
\begin{equation*}
	\widecheck{ c(a, k,M)} =^{(d)} a ^{H+\frac{1}{2}}I_{q} \left( \mathbbm{1}_{ (-\infty, 1) ^{q}\setminus (-M, 1) ^{q}}\int_{\mathbb{R}} dx \Psi (x) \int_{0} ^{x} f_{u}du\right).
\end{equation*}
Thus, by the isometry property (\ref{iso}) of the multiple stochastic integral $I_{q}$, we get
\begin{eqnarray}
	\mathbf{E} \Big[\widecheck{c (a, k,M)} ^{2}\Big] &=& a ^{2H+1} q! c_{q,H} ^{2}\int_{ \mathbbm{1}_{(-\infty, 1)^{q}\setminus (-M, 1) ^{q}}}  dy_{1}...dy_{q} \int_{\mathbb{R}} \int_{\mathbb{R}} dxdy \Psi (x) \Psi (y)\nonumber\\
	&& \int_{0} ^{x} dy \int_{0} ^{y} dv f_{u} (y_{1},...,y_{q}) f_{v}(y_{1},...,y_{q})\nonumber\\
	&\leq & C(q,H)a ^{2H+1} \int_{0} ^{1} \int_{0}^{1} dxdy \vert \Psi (x) \vert \cdot \vert \Psi (y) \vert  \int_{0} ^{x} du \int_{0} ^{y} dv \nonumber\\
	&&\times \left( \int_{\mathbb{R}} dy (u-y) _{+} ^{-\left( \frac{1}{2}+ \frac{1-H}{q}\right)}(v-y) _{+} ^{-\left( \frac{1}{2}+ \frac{1-H}{q}\right)}\right)^{q-1}\nonumber\\
	&& \left( \int_{-M} ^{1} dy (u-y) _{+} ^{-\left( \frac{1}{2}+ \frac{1-H}{q}\right)}(v-y) _{+} ^{-\left( \frac{1}{2}+ \frac{1-H}{q}\right)}\right)\nonumber\\
	&=& C(q,H)a ^{2H+1} \int_{0} ^{1} \int_{0}^{1} dxdy \vert \Psi (x) \vert \cdot \vert \Psi (y) \vert  \int_{0} ^{x} du \int_{0} ^{y} dv \vert u-v\vert ^{(2H-2) \frac{q-1}{q}}\nonumber\\
	&& \left( \int_{-M} ^{1} dy (u-y) _{+} ^{-\left( \frac{1}{2}+ \frac{1-H}{q}\right)}(v-y) _{+} ^{-\left( \frac{1}{2}+ \frac{1-H}{q}\right)}\right).\label{30n-2}
\end{eqnarray}
We have the following formula, obtained via the change of variables $z=\frac{ u\wedge v-y}{u\vee v-y}$,
\begin{eqnarray*}
	&& \int_{-M} ^{1} dy (u-y) _{+} ^{-\left( \frac{1}{2}+ \frac{1-H}{q}\right)}(v-y) _{+} ^{-\left( \frac{1}{2}+ \frac{1-H}{q}\right)}=\int_{-M} ^{u\wedge v} (u-y) ^{-\left( \frac{1}{2}+ \frac{1-H}{q}\right)}(v-y)  ^{-\left( \frac{1}{2}+ \frac{1-H}{q}\right)}dy \nonumber\\
	&=& \vert u-v\vert ^{\frac{2H-2}{q}}\int_{ 0} ^{ \frac{ u\wedge v+M}{u\vee v+M}}z ^{ \frac{H-1}{q} -\frac{1}{2}} (1-z) ^{\frac{2-2H}{q}-1} dz =\beta \left(  \frac{ u\wedge v+M}{u\vee v+M}, \frac{H-1}{q}+ \frac{1}{2}, \frac{2-2H}{q}\right),
\end{eqnarray*}
where $\beta (x, a, b)$ is the incomplete beta function given,  for $a,b>0$ and $x\in (0,1]$, by $ \beta (x, a, b)= \int_{0} ^{x} z ^{a-1} (1-z) ^{b-1}dz$. In particular, $\beta(1, a,b)= \beta (a,b)$  is the usual beta function. Let us also introduce the notation 
\begin{equation*}
	\bar{\beta} (x,a,b)= \beta (a,b)- \beta (x,a,b)= \int_{x} ^{1} z ^{a-1} (1-z) ^{b-1}dz
\end{equation*}
for $a,b>0$ and $ x \in (0,1).$  We have the following estimate 
\begin{eqnarray}
	&&	\bar{\beta} \left(  \frac{ u\wedge v+M}{u\vee v+M}, \frac{H-1}{q}+ \frac{1}{2}, \frac{2-2H}{q}\right)= \int_{ \frac{ u\wedge v+M}{u\vee v+M}} ^{1}z ^{\frac{H-1}{q}-\frac{1}{2}} (1-z) ^{\frac{2-2H}{q}-1}dz \nonumber \\
	&\leq & \int_{ \frac{M}{M+1}} ^{1} z ^{\frac{H-1}{q}-\frac{1}{2}} (1-z) ^{\frac{2-2H}{q}-1}dz\leq  \int_{ \frac{M}{M+1}} ^{1}(1-z) ^{\frac{2-2H}{q}-1}dz \nonumber \\
	&=&\frac{2-2H}{q}(M+1)^{\frac{2H-2}{q}}. \label{30n-1}
\end{eqnarray}
Therefore, by plugging (\ref{30n-1}) into (\ref{30n-2}), we get
\begin{eqnarray*}
	&&\mathbf{E} \Big[\widecheck{c(a, k,M) }^{2} \Big] 
	\\ &\leq& C(q,H) a ^{2H+1} (M+1)^{\frac{2H-2}{q}} \int_{0} ^{1} \int_{0}^{1} dxdy \vert \Psi (x) \vert \cdot \vert \Psi (y) \vert  \int_{0} ^{x} du \int_{0} ^{y} dv \vert u-v\vert ^{(2H-2) \frac{q-1}{q}}\\
	&\leq & C(q,H) a ^{2H+1} (M+1)^{\frac{2H-2}{q}} \int_{0} ^{1} \int_{0}^{1} dxdy \vert \Psi (x) \vert \cdot \vert \Psi (y) \vert  \int_{0} ^{x} du \int_{0} ^{y} dv \vert u-v\vert ^{2H-2}\\
	&=& C(q,H) a ^{2H+1} (M+1)^{\frac{2H-2}{q}} \int_{0} ^{1} \int_{0}^{1} dxdy \vert \Psi (x) \vert \cdot \vert \Psi (y) \vert  (x ^{2H}+ y ^{2H} -\vert x-y\vert ^{2H})\\
	&\leq &  C(q,H) a ^{2H+1} (M+1)^{\frac{2H-2}{q}}.
\end{eqnarray*}
\qed

\section{The modified wavelet variation} \label{sec:modifiedvar}

We will define and analyse the wavelet variation of the Hermite process based on some particular wavelet coefficients. We first introduce  these special  coefficients and we also study their main properties.

\subsection{The special wavelet coefficients}
We start with some notations. Consider two real numbers 
\begin{equation*}
	0<\gamma< \beta <1
\end{equation*}
and for every $N\geq 1$, let us set 
\begin{equation*}
	\mathcal{L}_{N}= \mathbb{N} \cap \left[ 1, \frac{ 2^{N}}{ 2 ^{[N ^{\beta}]}}\right]= \left[ 1, 2 ^{N-[N ^{\beta}]}\right]
\end{equation*}
and
\begin{equation*}
	\mathcal{L}_{N, \gamma} = \mathcal{L} _{N}\cap \left[ 1, 2 ^{ [N ^{\gamma}]}\right].
\end{equation*}
We have, for every $N\geq 1$, 
\begin{equation}\label{eq:cardiup}
	\mathcal{L}_{N, \gamma} \subset \mathcal{L} _{N} \mbox{ and } \vert \mathcal{L}_{N, \gamma}\vert\leq 2^{ [N ^{\gamma}]}
\end{equation}
and for $N$ large enough,
\begin{equation}\label{eq:cardi}
	\vert \mathcal{L}_{N, \gamma}\vert \geq 2^{ [N ^{\gamma}]}-1.
\end{equation}
For $ \ell \in \mathcal{L}_{N, \gamma}$ and $N\geq 1$, we set 
\begin{equation}
	\label{elnb}
	e_{\ell, N, \beta}= \frac{ \ell 2 ^{[N ^{\beta}]}}{2 ^{N}}.
\end{equation}

Now, we describe the wavelet coefficients on which we will focus in the rest of this work. Let $d\geq 1$. For $ \ell \in \mathcal{L}_{N, \gamma}$ and $ M=1,...,d$, we define
\begin{equation}
	\label{al}
	A_{M} (\ell, N)= c\left( \frac{1}{M2 ^{N}}, \ell M 2 ^{[N ^{\beta}]}\right). 
\end{equation}
From (\ref{coef1}), we express
\begin{eqnarray*}
		A_{M} (\ell, N)&=& \sqrt{ \frac{1}{ M  2 ^{N}}} I_{q} \left( 1^{\otimes q}_{ (-\infty, \frac{ \ell M 2 ^{[N ^{\beta}]}+M}{ M 2 ^{N}})}  \int_{\mathbb{R}} dx \Psi (x) \int_{ \frac{ \ell M  2 ^{[N ^{\beta}]}}{M ^{2 ^{N}}}} ^{{ \frac{ \ell M  2 ^{[N ^{\beta}]}+x}{M ^{2 ^{N}}}} } f_{u} du \right) \\&=&  \sqrt{ \frac{1}{ M  2 ^{N}}} I_{q} \left( 1 ^{\otimes q} _{ (-\infty, e_{ \ell, N, \beta}+ 2 ^{-N})}\int_{\mathbb{R} }dx \Psi (x) \int_{ e_{\ell, N, \beta}}^{e_{\ell, N, \beta}+xM ^{2 ^{-N}}}f_{u}du\right),
\end{eqnarray*}
where we used the notation (\ref{elnb}).  We have, by (\ref{2d-1}),
\begin{equation}
	\label{4d-3}
	\mathbf{E}\big[A_{M} (\ell, N)^{2}\big] =2 ^{-N (2H+1)} M ^{-(2H+1)}C_{\Psi }(H),
\end{equation}
with $C_{\Psi }(H) $ from (\ref{cpsi}). We will decompose this coefficient as
\begin{equation}\label{1d-2}
	A_{M} (\ell, N)= 	\widetilde{A_{M} (\ell, N)} + 	\widecheck{A_{M} (\ell, N)},
\end{equation}
where, with the notations (\ref{c1}) and (\ref{c2}),
\begin{eqnarray}
	\widetilde{A_{M} (\ell, N)} &=& \widetilde{c\left( \frac{1}{ M 2 ^{N}}, M \ell 2 ^{[N ^{\beta}]}, M ( 2 ^{[N ^{\beta}]}-1)\right)}\nonumber\\
	&=&  \sqrt{ \frac{1}{ M  2 ^{N}}} I_{q} \left( 1 ^{\otimes q} _{ (e_{ \ell-1, N, \beta}+ 2 ^{-N}, e_{ \ell, N, \beta}+ 2 ^{-N})}\int_{\mathbb{R} }dx \Psi (x) \int_{ e_{\ell, N, \beta}}^{e_{\ell, N, \beta}+xM ^{2 ^{-N}}}f_{u}du\right)\nonumber\\
	&=& I_{q}(\widetilde{g}_{M} (\ell, N)),\label{aln1}
\end{eqnarray}
with
\begin{eqnarray}
	\label{gnt}
	\widetilde{g}_{M} (\ell, N)(y_{1},...,y_{q})&=& \sqrt{ \frac{1}{ M  2 ^{N}}}  \mathbbm{1}^{\otimes q} _{ (e_{ \ell-1, N, \beta}+ 2 ^{-N}, e_{ \ell, N, \beta}+ 2 ^{-N})}(y_{1},...,y_{q})\\
	&&\int_{\mathbb{R} }dx \Psi (x) \int_{ e_{\ell, N, \beta}}^{e_{\ell, N, \beta}+xM ^{2 ^{-N}}}f_{u}(y_{1},...,y_{q})du\nonumber
\end{eqnarray}
and
\begin{eqnarray*}
&&	\widecheck{A_{M} (\ell, N)}\\
 &=& \widecheck{c\left( \frac{1}{ M 2 ^{N}}, M \ell 2 ^{[N ^{\beta}]}, M ( 2 ^{[N ^{\beta}]}-1)\right)}\\
	&=& \sqrt{ \frac{1}{ M  2 ^{N}}} I_{q} \left( \mathbbm{1}_{ (-\infty, e_{\ell, N, \beta}+2 ^{-N})^{q}\setminus   (e_{ \ell-1, N, \beta}+ 2 ^{-N}, e_{ \ell, N, \beta}+ 2 ^{-N})^{q}}\int_{\mathbb{R} }dx \Psi (x) \int_{ e_{\ell, N, \beta}}^{e_{\ell, N, \beta}+xM ^{2 ^{-N}}}f_{u}du\right)
\end{eqnarray*}

Before going further, let us emphasize the roles  for statistical inference of the various parameters introduced in this section. The estimator $\widehat{H}_N$ for the Hurst parameter defined in equation \eqref{est1} below relies on the wavelet variations \eqref{vnl} in which the corresponding wavelet coefficients are precisely chosen. First, we need to consider $d$ scales, that is the parameter $a$ in \eqref{coef}, to perform a $\log$-regression. To define $\widehat{H}_N$, we take the scales $(M2^N)^{-1}$ for $M=1,\dots,d$. Then, at such a scale $(M2^N)^{-1}$ we only work with wavelet coefficients at position, that is the parameter $k$ in \eqref{coef}, which are multiple of  $M 2 ^{[N ^{\beta}]}$. Thus, we exactly get the definition of the coefficients $A_M(\ell,N)$ in \eqref{al}. This choice of wavelet coefficients has the double advantage to offer, as already stated, some independence properties for the ``tilde part'', see also Lemma \ref{ll2} below, while the ``check counterpart'' are negligible, for the $L^2(\Omega)$-norm. Nevertheless, if all coefficients $A_M(\ell,N)$, with $\ell \in \mathcal{L}_{N}$, were chosen in \eqref{vnl}, then the proofs of the Central Limit Theorems (Theorem \ref{tt1} and \ref{thm:2} below) would unfortunately not hold. To overcome this situation, we just have to select a proportion of the coefficients $A_M(\ell,N)$, that is the one for which $\ell \in \mathcal{L}_{N, \gamma}$. The simple fact that $\gamma <  \beta$ allows to deduce precious convergence to $0$ in $L^1(\Omega)$, see for instance Proposition \ref{pp1} below. Note that increasing $\beta$ allows to also increase $\gamma$ and thus improve the rate of convergence to the normal distribution, see the bound for the Wasserstein distance in Theorems \ref{tt1} and \ref{thm:2}, but, on the other side, reduce the number of data used in the statistical inference, as a consequence of the definition of $\mathcal{L}_{N}$ and $\mathcal{L}_{N, \gamma}$. For this reason, it is a difficult question to know what should be the more ``practicable'' value for the parameter $\beta$.

Let us come back to the main properties of the special wavelet coefficients introduced here over. First, the following facts will play an important role for the proofs of our main results. 

\begin{lemma}\label{ll2}
\mbox {}
\begin{enumerate}
	\item 
	For every $N\geq 1$, the random vectors $\left( ( \widetilde{A_{1}(\ell, N)}, \ldots, \widetilde{A_{d}(\ell, N)}), \ell \in \mathcal{L}_{N, \gamma}\right)$ are mutually independent.
	
	\item For every $N\geq 1$ and $M=1,..., d$, the random variables $( \widetilde{A_{M}(\ell, N)}, \ell \in \mathcal{L}_{N, \gamma})$ are identically distributed. 
	
		\item For every $N\geq 1$ and $M=1,..., d$, the random variables $( \widecheck{A_{M}(\ell, N)}, \ell \in \mathcal{L}_{N, \gamma})$ are identically distributed. 
	
\end{enumerate}
\end{lemma} 
\noindent {\bf Proof: } To prove the first point, we notice that for every $\ell, j \in \mathcal{L}_{N, \gamma}$ with $\ell\not=j$ and for every $M_{1}, M_{2}=1,...,d$, we clearly have 
\begin{equation*}
	\widetilde{g_{M_{1}} (\ell, N)}\otimes _{1} \widetilde{g_{M_{2}} (j, N)}=0
\end{equation*}
almost everywhere on $ \mathbb{R} ^{2q-2}$ (the notation $\otimes_{1}$ stands for the contraction of order one, see (\ref{contra}) in the Appendix).  This is because the intervals that appear in the expression of the kernels $\widetilde{g_{M_{1}} (\ell, N)}$ and $\widetilde{g_{M_{2}} (j, N)}$ are disjoint. The independence follows from the Ust\'unel-Zakai criterion (see e.g. \cite[Proposition 1]{UZ}). 

Concerning the second point, for any integer $h$ such that $\ell+h\in \mathcal{L}_{N, \gamma}$, we can write 
\begin{eqnarray*}
	\widetilde{	A_{M}(\ell +h, N)}&=& d(H) \sqrt{ \frac{1}{M2 ^{N}}} \int_{\mathbb{R} ^{q}} dB(y_{1})...dB(y_{q})\\
	&& \mathbbm{1}_{ ( e_{\ell+h-1, N, \beta}+ 2 ^{-N}, e_{\ell+h, N, \beta}+ 2 ^{-N})}^{\otimes q}(y_{1},..., y_{q})\\
	&&
	\int_{\mathbb{R}} dx \Psi (x) \int _{ e_{\ell+h, N, \beta}}^{e_{\ell+h, N, \beta}+ x M2 ^{-N}}(u-y_{1})_{+} ^{-\left( \frac{1}{2}+\frac{1-H}{q}\right)}\ldots (u-y_{q})_{+} ^{-\left( \frac{1}{2}+\frac{1-H}{q}\right)}du 
\end{eqnarray*}
and with the successive changes of variables $\tilde{u}=u-e_{h, N, \beta}$ and $ \tilde{y}_{j}= y_{j}- e_{h, N, \beta}$ for $ j=1,...,q$, we arrive at 
\begin{eqnarray*}
&&	\widetilde{	A_{M}(\ell +h, N)}\\
&=& d(H) \sqrt{ \frac{1}{M2 ^{N}}} \int_{\mathbb{R} ^{q}} dB(y_{1})...dB(y_{q})\\
	&& \mathbbm{1}_{ ( e_{\ell+h-1, N, \beta}+ 2 ^{-N}, e_{\ell+h, N, \beta}+ 2 ^{-N})}^{\otimes q}(y_{1},..., y_{q})\\
	&&\int_{\mathbb{R}} dx \Psi (x) \int _{ e_{\ell, N, \beta}}^{e_{\ell, N, \beta}+ x M2 ^{-N}}(u-y_{1}+ e_{h, N, \beta})_{+} ^{-\left( \frac{1}{2}+\frac{1-H}{q}\right)}\ldots (u-y_{q}+ e_{h, N, \beta})_{+} ^{-\left( \frac{1}{2}+\frac{1-H}{q}\right)}du \\
	&=& d(H)  \sqrt{ \frac{1}{M2 ^{N}}} \int_{\mathbb{R} ^{q}} dB(y_{1}+ e_{h, N, \beta})...dB(y_{q}+ e_{h, N, \beta})\\
	&& \mathbbm{1}_{ ( e_{\ell-1, N, \beta}+ 2 ^{-N}, e_{\ell, N, \beta}+ 2 ^{-N})}^{\otimes q}(y_{1},..., y_{q})\\
	&&\int_{\mathbb{R}} dx \Psi (x) \int _{ e_{\ell, N, \beta}}^{e_{\ell, N, \beta}+ x M2 ^{-N}}(u-y_{1})_{+} ^{-\left( \frac{1}{2}+\frac{1-H}{q}\right)}\ldots (u-y_{q})_{+} ^{-\left( \frac{1}{2}+\frac{1-H}{q}\right)}du.
\end{eqnarray*}
Since the Brownian motion has stationary increments, we have 
\begin{eqnarray*}
	\widetilde{	A_{M}(\ell +h, N)}&=^{(d)}& d(H)  \sqrt{ \frac{1}{M2 ^{N}}} \int_{\mathbb{R} ^{q}} dB(y_{1})...dB(y_{q})\\
	&& \mathbbm{1}_{ ( e_{\ell-1, N, \beta}+ 2 ^{-N}, e_{\ell, N, \beta}+ 2 ^{-N})}^{\otimes q}(y_{1},..., y_{q})\\
	&&\int_{\mathbb{R}} dx \Psi (x) \int _{ e_{\ell, N, \beta}}^{e_{\ell, N, \beta}+ x M2 ^{-N}}(u-y_{1})_{+} ^{-\left( \frac{1}{2}+\frac{1-H}{q}\right)}\ldots (u-y_{q})_{+} ^{-\left( \frac{1}{2}+\frac{1-H}{q}\right)}du \\
	&=& 	\widetilde{	A_{M}(\ell , N)}.
\end{eqnarray*}
The third point follows immediately, using the same arguments as above. \qed

The next result is then a direct consequence of Lemma \ref{ll1}.

\begin{lemma}\label{ll3}
	For every $N\geq 1, M=1,...,d$ and for every $\ell \in \mathcal{L}_{N,\gamma}$, we have
	\begin{equation*}
		\mathbf{E}\Big[ \widecheck{A_{M}(\ell, N)}^{2} \Big] \leq C(q,H, M) 2 ^{-N(2H+1)} 2 ^{ [N ^{\beta}]\frac{2H-2}{q}}.
	\end{equation*}

\end{lemma}

We will need some auxiliary results concerning the behavior of the second and fourth moments of the random variable (\ref{aln1}).

\begin{lemma}\label{ll4} For $N$ sufficiently large and for all $\ell \in \mathcal{L}_{N, \gamma}, M=1,...,d$, we have 
	\begin{equation}
		\label{1d-1}
		\left| 2 ^{N (2H+1)} \mathbf{E} \Big[\widetilde{	A_{M}(\ell , N)}^{2} \Big]- M ^{-(2H+1)}C_{\Psi} (H) \right| \leq C2 ^{ [N ^{\beta}]\frac{2H-2}{q}}.
	\end{equation}
In particular,
\begin{equation*}
	2 ^{N (2H+1)} \mathbf{E} \Big[\widetilde{	A_{M}(\ell , N)}^{2} \Big] \underset{{N \to \infty}}{\longrightarrow} M ^{-(2H+1)}C_{\Psi} (H).
\end{equation*}
Also, for $N$ sufficiently large and for $\ell \in \mathcal{L}_{N, \gamma}, M=1,...,d$, we have 
\begin{equation}
	\label{1d-4}
		\left| 2 ^{2N (2H+1)} \mathbf{E} \Big[\widetilde{	A_{M}(\ell , N)}^{4} \Big] - M ^{-(4H+2)} \mathbf{E} \big[ c(1,0)^{4} \big]\right| \leq C2 ^{ [N ^{\beta}]\frac{2H-2}{q}},
\end{equation}
which implies
\begin{equation*}
	2 ^{2N (2H+1)} \mathbf{E} \Big[\widetilde{	A_{M}(\ell , N)}^{4}\Big] {{N \to \infty}} M ^{-(4H+2)}\mathbf{E} \big[c(1,0) ^{4}\big].
\end{equation*}
\end{lemma}
\noindent {\bf Proof: } Since $\widetilde{	A_{M}(\ell , N)}$ and $ \widecheck{A_{M}(\ell, N)}$ are independent random variables, we have from (\ref{1d-2}),
\begin{equation}\label{4d-1}
	\mathbf{E} \Big[A_{M} (\ell, N) ^{2}\Big]= \mathbf{E} \Big[\widetilde{	A_{M}(\ell , N)}^{2}\Big] +\mathbf{E} \Big[\widecheck{	A_{M}(\ell , N)}^{2}\Big].
\end{equation}
Let us note that, the above inequality combined with equality \eqref{4d-3} and Lemme \ref{ll3} entail
\begin{equation}\label{4d-2}
2 ^{N (2H+1)} \mathbf{E} \Big[\widetilde{	A_{M}(\ell , N)}^{2}\Big] \leq C.
\end{equation}

The relation (\ref{4d-1}) gives 
\begin{eqnarray*}
	2 ^{N (2H+1)} \mathbf{E} \Big[\widetilde{	A_{M}(\ell , N)}^{2} \Big]&=&2 ^{N (2H+1)} 	\mathbf{E}\Big[ A_{M} (\ell, N) ^{2}\Big]- 2 ^{N (2H+1)}\mathbf{E} \Big[\widecheck{	A_{M}(\ell , N)}^{2} \Big]\\
&=& M ^{-(2H+1)}C_{\Psi} (H)- 2 ^{N (2H+1)}\mathbf{E} \Big[\widecheck{	A_{M}(\ell , N)}^{2}\Big],
\end{eqnarray*}
and the inequality (\ref{1d-1}) is then a direct consequence of Lemma \ref{ll3}. 

To prove (\ref{1d-4}), we write
\begin{eqnarray*}
&&	\mathbf{E} 	\Big[A_{M}(\ell , N)^{4}\Big]\\
&=& \mathbf{E} \Big[\widetilde{	A_{M}(\ell , N)}^{4}\Big]+ \mathbf{E} \Big[\widecheck{	A_{M}(\ell , N)}^{4}\Big]+ 6\mathbf{E} \Big[\widetilde{	A_{M}(\ell , N)}^{2}\widecheck{	A_{M}(\ell , N)}^{2}\Big]\\
	&&+ 4\mathbf{E} \Big[\widetilde{	A_{M}(\ell , N)}^{3}\widecheck{	A_{M}(\ell , N)}\Big]+ 4\mathbf{E} \Big[\widetilde{	A_{M}(\ell , N)}\widecheck{	A_{M}(\ell , N)}^{3}\Big]\\
	&=& \mathbf{E} \Big[\widetilde{	A_{M}(\ell , N)}^{4}\Big]+ \mathbf{E} \Big[\widecheck{	A_{M}(\ell , N)}^{4}\Big]+ 6\mathbf{E} \Big[\widetilde{	A_{M}(\ell , N)}^{2}\widecheck{	A_{M}(\ell , N)}^{2}\Big]+4\mathbf{E} \Big[\widetilde{	A_{M}(\ell , N)}\widecheck{	A_{M}(\ell , N)}^{3}\Big],
\end{eqnarray*}
because, by \cite[Lemma 2]{AT}, we have $\mathbf{E} \Big[ \widetilde{	A_{M}(\ell , N)}^{3}\widecheck{	A_{M}(\ell , N)}\Big]=0$. Thus, we get, using  the equality (\ref{2d-1}),
\begin{eqnarray*}
	&&2 ^{2N (2H+1)} \mathbf{E} \Big[\widetilde{	A_{M}(\ell , N)}^{4}\Big]\\
	& = &2 ^{2N (2H+1)}	\mathbf{E} \Big[	A_{M}(\ell , N)^{4}\Big]-2^{2N (2H+1)}  \mathbf{E} \Big[\widecheck{	A_{M}(\ell , N)}^{4}\Big]\\
	&&- 6 \times 2 ^{2N (2H+1)}\mathbf{E} \Big[\widetilde{	A_{M}(\ell , N)}^{2}\widecheck{	A_{M}(\ell , N)}^{2} \Big]-4 \times 2 ^{2N (2H+1)} \mathbf{E} \Big[\widetilde{	A_{M}(\ell , N)}\widecheck{	A_{M}(\ell , N)}^{3}\Big]\\
	&=&M ^{-(4H+2)}\mathbf{E} \Big[c(1,0)^{4}\Big]-2^{2N (2H+1)}  \mathbf{E} \Big[\widecheck{	A_{M}(\ell , N)}^{4}\Big]\\
	&&- 6 \times 2 ^{2N (2H+1)}\mathbf{E} \Big[\widetilde{	A_{M}(\ell , N)}^{2}\widecheck{	A_{M}(\ell , N)}^{2}\Big]-4  \times 2 ^{2N (2H+1)} \mathbf{E} \Big[\widetilde{	A_{M}(\ell , N)} \widecheck{	A_{M}(\ell , N)}^{3}\Big].
\end{eqnarray*}
The hypercontractivity property \eqref{hyper}, the inequality (\ref{4d-2})  and Lemma \ref{ll3} imply
\begin{equation*}
	2^{2N (2H+1)}  \mathbf{E} \Big[\widecheck{	A_{M}(\ell , N)}^{4}\Big]\leq C 	2^{2N (2H+1)}  \left(\mathbf{E} \Big[\widecheck{	A_{M}(\ell , N)}^{2}\Big]\right) ^{2} \leq C 2 ^{[ N ^{\beta}]\frac{ 4H-4}{q}}
\end{equation*}
 while, by, the Cauchy-Schwarz's inequality, \eqref{hyper}, Lemma \ref{ll3} and (\ref{4d-2}), we have
\begin{eqnarray*}
\mathbf{E} \Big[\widetilde{	A_{M}(\ell , N)}^{2}\widecheck{	A_{M}(\ell , N)}^{2}\Big] & \leq &  \left(\mathbf{E} \Big[\widetilde{	A_{M}(\ell , N)}^{4} \Big] \right)^\frac{1}{2} \left( \mathbf{E}\Big[\widecheck{	A_{M}(\ell , N)}^{4}\Big] \right)^\frac{1}{2} \\
& \leq& C  \mathbf{E} \Big[\widetilde{	A_{M}(\ell , N)}^{2} \Big]  \mathbf{E}\Big[\widecheck{	A_{M}(\ell , N)}^{2}\Big] \\
&\leq& C 2 ^{-2N (2H+1)} 2 ^{[ N ^{\beta}]\frac{ 2H-2}{q}}
\end{eqnarray*}
and
\begin{eqnarray*}
	 \mathbf{E} \Big[\widetilde{	A_{M}(\ell , N)}\widecheck{	A_{M}(\ell , N)}^{3}\Big]&\leq &C \left( \mathbf{E} \Big[\widetilde{	A_{M}(\ell , N)}^{2}\Big] \right) ^{\frac{1}{2}}\left( \mathbf{E} \Big[\widecheck{	A_{M}(\ell , N)}^{6}\Big]\right) ^{\frac{1}{2}}\\
	 &\leq  &C \left( \mathbf{E} \Big[\widetilde{	A_{M}(\ell , N)}^{2}\Big] \right) ^{\frac{1}{2}}\left( \mathbf{E} \Big[\widecheck{	A_{M}(\ell , N)}^{2}\Big]\right) ^{\frac{3}{2}}\\
	 &\leq& C 2 ^{-2N (2H+1)} 2 ^{[ N ^{\beta}]\frac{ 3H-3}{q}}.
\end{eqnarray*}
The above three bounds lead to inequality (\ref{1d-4}). \qed 

\begin{lemma}\label{l5}
	Let $M_{1}, M_{2} =1,...,d$. Then for $N$ sufficiently large, we have 
	\begin{equation*}
	\left| 2 ^{2N(2H+1)}\mathbf{E} \Big[\widetilde{A_{M_{1}}(\ell, N) }^{2} \widetilde{A_{M_{2}}(\ell, N) }^{2}\Big]- \mathbf{E} \Big[c(M_{1}^{-1}, 0)^{2} c(M_{2}^{-1}, 0) ^{2}\Big]\right| \leq C 2 ^{[ N ^{\beta}]\frac{ 2H-2}{q}}.
	\end{equation*}
\end{lemma}
\noindent { \bf Proof: } By the scaling property and the stationarity of the increments of the Hermite process we can show that for every $M_{1}, M_{2}=1,...,d$, 
\begin{equation*}
	( A_{M_{1}}(\ell, N), A_{M_{2}}(\ell, N))=^{(d)} 2 ^{-N(H+\frac{1}{2})}( c(M_{1}^{-1},0), c(M_{2}^{-1},0)).
\end{equation*}
This implies
\begin{equation*}
	2 ^{2N(2H+1)}\mathbf{E} \Big[A_{M_{1}}(\ell, N)^{2} A_{M_{2}}(\ell, N) ^{2}\Big]= \mathbf{E} \Big[c(M_{1}^{-1}, 0)^{2} c(M_{2}^{-1}, 0) ^{2}\Big].
\end{equation*}
By using the decomposition (\ref{1d-1}) for $ A_{M_{1}}(\ell, N)$ and $ A_{M_{2}}(\ell, N)$ and by applying the estimate  in Lemma \ref{ll3} for the negligible parts $\widecheck{ A_{M_{1}}(\ell, N)}$ and $\widecheck{ A_{M_{2}}(\ell, N)}$, we can conclude. \qed

\subsection{The modified wavelet variation and its negligible part}
Let us now introduce the main object of this work, namely the modified wavelet variation of the Hermite process. This wavelet variation is defined by using, not all the wavelet coefficients of the Hermite process, but only the special coefficient given by (\ref{al}). More exactly, for $N\geq 1, M=1,...,d$, we set

\begin{equation}
	\label{vnl}
	V_{N, M}=\frac{1}{\sqrt{\vert \mathcal{L}_{N, \gamma}\vert}}\sum_{\ell \in \mathcal{L}_{N, \gamma}}\left[ \frac{ A_{M} (\ell, N)^{2} }{ \mathbf{E} \big[A_{M}(\ell, N) ^{2}\big]}-1\right].
\end{equation}
From (\ref{4d-3}), we can also write 
\begin{equation}
	\label{vnl2}
	V_{N, M}=\frac{1}{\sqrt{\vert \mathcal{L}_{N, \gamma}\vert }} \frac{ 2 ^{N (2H+1)} M ^{2H+1}}{C_{\Psi}(H)}\sum_{ \ell \in \mathcal{L}_{N, \gamma} }\left( A_{M}(\ell, N) ^{2}- \mathbf{E}\big[A_{M}(\ell, N) ^{2}\big]\right). 
\end{equation}
The purpose is to find the limit behavior in distribution, as $ N \to \infty$, of the sequence $ (V_{N, M}, N\geq 1)$. To this end, we decompose $ V_{N, M}$ as follows
\begin{equation}
	\label{deco}
	V_{N, M}= V_{N, M, 1}+ V_{N, M, 2}+ V_{N,M, 3},
\end{equation}
where, for all $N\geq 1$ and $M=1,...,d$, we set
\begin{equation}
	\label{v1}
	V_{N, M, 1}=\frac{1}{\sqrt{\vert \mathcal{L}_{N, \gamma}\vert }} \frac{ 2 ^{N (2H+1)} M ^{2H+1}}{C_{\Psi}(H)}\sum_{ \ell \in \mathcal{L}_{N, \gamma} }\left( \widetilde{A_{M}(\ell, N)} ^{2}- \mathbf{E}\Big[\widetilde{A_{M}(\ell, N)} ^{2} \Big]\right),
\end{equation}

\begin{equation}
	\label{v2}
		V_{N, M, 2}=\frac{1}{\sqrt{\vert \mathcal{L}_{N, \gamma}\vert }} \frac{ 2 ^{N (2H+1)} M ^{2H+1}}{C_{\Psi}(H)}\sum_{ \ell \in \mathcal{L}_{N, \gamma} }\left( \widecheck{A_{M}(\ell, N)} ^{2}- \mathbf{E}\Big[\widecheck{A_{M}(\ell, N)} ^{2}\Big]\right),
\end{equation}
and
\begin{equation}
	\label{v3}
	V_{N, M, 3}= 2\frac{1}{\sqrt{\vert \mathcal{L}_{N, \gamma}\vert }} \frac{ 2 ^{N (2H+1)} M ^{2H+1}}{C_{\Psi}(H)}\sum_{ \ell \in \mathcal{L}_{N, \gamma} }\widetilde{A_{M}(\ell, N)}\widecheck{A_{M}(\ell, N)}.
\end{equation}
In a first step, we will show that the summands denoted by $ V_{N, M, 2}$ and $ V_{N, M, 3} $ above are negligible, and consequently the behavior of the wavelet variation  (\ref{vnl}) is given by the sequence $(V_{N, M, 1}, N\geq 1)$. 

\begin{prop}\label{pp1}
	Let $ V_{N, M,2}, V_{N, M, 3}$ be given by (\ref{v2}), (\ref{v3}), respectively. We have \begin{equation*}
		\mathbf{E} \big[\vert V_{N, M, 2} \vert \big] \leq C(q,H, M) 2 ^{ \frac{ N ^{\gamma}}{2}} 2 ^{ [N ^{\beta}]\frac{2H-2}{q}}
	\end{equation*}
	and
	\begin{equation*}
		\mathbf{E} \big[\vert V_{N, M, 3}\vert \big]  \leq C(q,H, M) 2 ^{ \frac{ N ^{\gamma}}{2}} 2 ^{ [N ^{\beta}]\frac{H-1}{q}}.
	\end{equation*}
In particular, the sequences $(V_{N,M,2}, N\geq 1)$ and $(V_{N,M,3}, N\geq 1)$ converge to zero in $ L ^{1}(\Omega)$, as $ N \to \infty$.
\end{prop}
\noindent { \bf Proof: } For $V_{N,M,2}$ we have the following  estimates, by (\ref{v2}) and  Lemma \ref{ll3},
\begin{eqnarray*}
	\mathbf{E} \big[\vert V_{N,M,2}\vert \big] &\leq & C(q,H, M) \frac{ 2 ^{N (2H+1)}}{\sqrt{ \vert \mathcal{L}_{N, \gamma}\vert }}\sum_{ \ell \in \mathcal{L}_{N, \gamma}}\mathbf{E}\Big[\widecheck{A_{M}(\ell, N) }^{2}\Big]\\
	&\leq & C(q,H, M) \sqrt{ \vert \mathcal{L}_{N, \gamma}\vert } 2 ^{-N (2H+1)}2 ^{ [N ^{\beta}]\frac{2H-2}{q}}\\
	&\leq & C(q,H, M)2 ^{ \frac{ N ^{\gamma}}{2}}2 ^{ [N ^{\beta}]\frac{2H-2}{q}}.
\end{eqnarray*}
Since $\gamma <\beta$ and $H<1$, we get the converge to zero in $ L ^{1}(\Omega)$, as $N\to \infty$,  of $V_{N,M,2}$.

Concerning the summand denoted by $ V_{N,M, 3}$, we have by (\ref{v3}) and Cauchy-Schwarz's inequality,

\begin{eqnarray*}
	\mathbf{E} \big[\vert V_{N,M,3}\vert \big] &\leq & C(q,H, M) \frac{ 2 ^{N (2H+1)}}{\sqrt{ \vert \mathcal{L}_{N, \gamma}\vert }}\sum_{ \ell \in \mathcal{L}_{N, \gamma}}\left( \mathbf{E} \Big[\widetilde{A_{M}(\ell, N) }^{2}\Big]\right) ^{\frac{1}{2}} \left( \mathbf{E}\Big[\widecheck{A_{M}(\ell, N) }^{2}\Big]\right) ^{\frac{1}{2}}\\
	&\leq & C(q,H, M) \frac{ 2 ^{N (2H+1)}}{\sqrt{ \vert \mathcal{L}_{N, \gamma}\vert }}\sum_{ \ell \in \mathcal{L}_{N, \gamma}} \sqrt{ 2 ^{N (2H+1)}} \sqrt{ 2 ^{N (2H+1)}2 ^{ [N ^{\beta}]\frac{2H-2}{q}}}\\
	&\leq & C(q,H, M)  2 ^{ \frac{ N ^{\gamma}}{2}} 2 ^{ [N ^{\beta}]\frac{H-1}{q}}.
\end{eqnarray*}
Again, the assumptions $\gamma <\beta$ and $H<1$ imply that $V_{N,M,3}$ converges to zero in $ L ^{1}(\Omega)$ as $ N \to \infty$. \qed

\subsection{The Central Limit Theorem for the modified wavelet variation}

In this section, we prove that the $d$-dimensional random vector $ ( V_{N, M}, M=1,...,d)$ satisfies a (multidimensional) Central Limit Theorem. On this purpose, we first focus on the sequence $(V_{N, M, 1}, N\geq 1)$ defined by (\ref{v1}).  

Let us start by evaluating the behavior of $\mathbf{E}\big[ V_{N, M_{1}, 1} V_{N, M_{2}, 1}\big]$, with $M_1,M_2=1,...,d$, as $N\to \infty$. 

\begin{prop}\label{pp2}
	Let $ V_{N, M, 1}$ be given by (\ref{v1}). Then for $N$ large enough, and for every $M=1,...,d$ , we have
	\begin{equation}\label{6d-4}
		\left| \mathbf{E}\big[ V_{N, M, 1} ^{2}\big]-  \left(\frac{ \mathbf{E}\big[ c(1,0) ^{4}\big]}{ C_{\Psi}(H) ^{2}}- 1 \right)\right|  \leq C2 ^{[ N ^{\beta}]\frac{ 2H-2}{q}}.
	\end{equation}
If $M_{1}, M_{2}=1,...,d$ are such that $M_{1}\not=M_{2}$, then for $N$ large 
\begin{equation}\label{8d-1}
	\left| \mathbf{E}\big[ V_{N, M_{1}, 1} V_{N, M_{2}, 1}\big]- \left(\frac{(M_{1}M_{2})^{2H+1}\mathbf{E}\big[ c(M_{1}^{-1},0) ^{2}c(M_{2}^{-1}, 0)^{2}\big]}{C_{\Psi}(H) ^{2}} -1 \right)\right| \leq  C2 ^{[ N ^{\beta}]\frac{ 2H-2}{q}}.
\end{equation}
\end{prop}
\noindent {\bf Proof: } By (\ref{v1}), we have, for $M=1,...,d$ and $N\geq 1$, 
\begin{eqnarray*}
	 V_{N, M, 1} ^{2}&=& \frac{1}{ \vert \mathcal{L}_{N, \gamma}\vert }\frac{ 2 ^{2N (2H+1)M ^{4H+2}}}{C_{\Psi} (H)^{2}}\\
	&&\sum _{\ell, j\in \mathcal{L} _{N, \gamma}}\mathbf{E}\left[ \left( \widetilde{A_{M}(\ell, N)} ^{2}- \mathbf{E}\big[\widetilde{A_{M}(\ell, N) }^{2}\big]\right)\left( \widetilde{A_{M}(j, N)} ^{2}- \mathbf{E}\big[\widetilde{A_{M}(j, N)} ^{2}\big]\right)\right]\\
	&=&  \frac{1}{ \vert \mathcal{L}_{N, \gamma}\vert }\frac{ 2 ^{2N (2H+1)M ^{4H+2}}}{C_{\Psi} (H)^{2}}\sum _{\ell\in \mathcal{L} _{N, \gamma}}\mathbf{E}\left[\left( \widetilde{A_{M}(\ell, N) }^{2}- \mathbf{E}\big[\widetilde{A_{M}(\ell, N)} ^{2}\big]\right)^{2}\right],
\end{eqnarray*}
since $ \widetilde{A_{M} (\ell, N)}$ and $ \widetilde{A_{M}(j, N)}$ are independent when $\ell\not=j$ (see the first point of Lemma \ref{ll2}). Next, since, by the second point of Lemma \ref{ll2}, $ (\widetilde{A_{M}(\ell, N)}, \ell \in \mathcal{L}_{N, \gamma})$ are identically distributed,  we can write 
\begin{eqnarray*}
	\mathbf{E}\big[ V_{N, M, 1} ^{2}\big]&=& \frac{ 2 ^{2N (2H+1)M ^{4H+2}}}{C_{\Psi} (H)^{2}}\mathbf{E}\left[ \left( \widetilde{A_{M}(\ell_{0}, N) }^{2}- \mathbf{E}\Big[ \widetilde{A_{M}(\ell_{0}, N)} ^{2}\Big]\right)^{2}\right]\\
	&=& \frac{ 2 ^{2N (2H+1)M ^{4H+2}}}{C_{\Psi} (H)^{2}}\left( \mathbf{E} \Big[\widetilde{A_{M}(\ell_{0}, N)}^{4}\Big] - \left(\mathbf{E}  \Big[\widetilde{A_{M}(\ell_{0}, N)}^{2}\Big]\right) ^{2} \right).
\end{eqnarray*}
for some $\ell_{0} \in \mathcal{L}_{N, \gamma}$.  Therefore, we obtain
\begin{eqnarray}
	&&\mathbf{E}[ V_{N, M, 1} ^{2}]-  \left(\frac{ \mathbf{E}\big[ c(1,0) ^{4}\big]}{ C_{\Psi}(H) ^{2}}- 1 \right)\nonumber\\
	&=&\frac{1}{ C_{\Psi} (H)^{2}} M ^{4H+2}\left[  \left(  2 ^{2N(2H+1)}\mathbf{E}  \Big[\widetilde{A_{M}(\ell_{0}, N)}^{4}\Big]- M ^{-(4H+2)}\mathbf{E}[ c(1,0) ^{4}]\right) \right.\nonumber \\
	&&\left. - \left( 2 ^{2N(2H+1)}\left( \mathbf{E} \Big[\widetilde{A_{M}(\ell_{0}, N)}^{2}\Big]\right)^2 - M ^{-(4H+2)}C_\Psi(H)^{2} \right) \right].\label{4d-4}
\end{eqnarray}
On one hand, by (\ref{1d-4}), we have
\begin{equation}
	\label{4d-5}
\left| 	 2 ^{2N(2H+1)}\mathbf{E}  \Big[\widetilde{A_{M}(\ell_{0}, N)}^{4}\Big]- M ^{-(4H+2)}\mathbf{E} \big[c(1,0) ^{4}\big]\right| \leq C2 ^{[ N ^{\beta}]\frac{ 2H-2}{q}}
\end{equation}
while, on the other hand,  by Lemma \ref{ll4}, we get
\begin{eqnarray}
&&\left| 	2 ^{2N(2H+1)}\left( \mathbf{E}  \Big[\widetilde{A_{M}(\ell_{0}, N)}^{2}\Big]\right)^2- M ^{-(4H+2)}C_\Psi(H)^{2}\right| \nonumber\\
&& = \left|  2 ^{N(2H+1)}\mathbf{E}  \Big[\widetilde{A_{M}(\ell_{0}, N)}^{2}\Big]- M ^{-(2H+1)}C_\Psi(H)\right| \left|  2 ^{N(2H+1)}\mathbf{E} \Big[ \widetilde{A_{M}(\ell_{0}, N)}^{2}\Big]+ M ^{-(2H+1)}C_\Psi(H)\right|\nonumber \\
&&\leq C 2 ^{[ N ^{\beta}]\frac{ 2H-2}{q}}.\label{4d-6}
\end{eqnarray}
By plugging (\ref{4d-5}) and (\ref{4d-6}) into (\ref{4d-4}), we obtain the inequality (\ref{6d-4}). 

If $M_{1}\not=M_{2}$, by (\ref{vnl2}) and the independence property in the first point of Lemma \ref{ll2}, we express
\begin{eqnarray*}
&&	\mathbf{E}\big[V_{N, M_{1}, 1} V_{N, M_{2}, 1}\big]\\
&=&  \frac{1}{ \vert \mathcal{L}_{N, \gamma}\vert }\frac{ 2 ^{2N (2H+1)M_{1} ^{2H+1}M_{2} ^{2H+1}}}{C_{\Psi} (H)^{2}}\sum _{\ell\in \mathcal{L} _{N, \gamma}}\\
	&&\mathbf{E}\left[\left( \widetilde{A_{M_{1}}(\ell, N) }^{2}- \mathbf{E}\Big[\widetilde{A_{M_{1}}(\ell, N)} ^{2}\Big]\right)\left( \widetilde{A_{M_{2}}(\ell, N) }^{2}- \mathbf{E}\Big[\widetilde{A_{M_{2}}(\ell, N)} ^{2}\Big]\right)\right]\\
	&=&\frac{1}{ \vert \mathcal{L}_{N, \gamma}\vert }\frac{ 2 ^{2N (2H+1)M_{1} ^{2H+1}M_{2} ^{2H+1}}}{C_{\Psi} (H)^{2}}\sum _{\ell\in \mathcal{L} _{N, \gamma}}\\
&&	\left( \mathbf{E}\Big[\widetilde{A_{M_{1}}(\ell, N) }^{2}\widetilde{A_{M_{2}}(\ell, N) }^{2}\Big]- \mathbf{E}\Big[\widetilde{A_{M_{1}}(\ell, N) }^{2}\Big]\mathbf{E}\Big[\widetilde{A_{M_{2}}(\ell, N) }^{2}\Big]\right)\\
&=&\frac{ 2 ^{2N (2H+1)M_{1} ^{2H+1}M_{2} ^{2H+1}}}{C_{\Psi} (H)^{2}}	\left( \mathbf{E}\Big[\widetilde{A_{M_{1}}(\ell_0, N) }^{2}\widetilde{A_{M_{2}}(\ell_0, N) }^{2}\Big]- \mathbf{E}\Big[\widetilde{A_{M_{1}}(\ell_0, N) }^{2}\Big]\mathbf{E}\Big[\widetilde{A_{M_{2}}(\ell_0, N) }^{2} \Big]\right)
\end{eqnarray*}
for some $\ell_0 \in \mathcal{L}_{N, \gamma}$. This time, we write

\begin{eqnarray}
&&	\mathbf{E}\big[V_{N, M_{1}, 1} V_{N, M_{2}, 1}\big]-\left(\frac{(M_{1}M_{2})^{2H+1}\mathbf{E}\big[ c(M_{1}^{-1},0) ^{2}c(M_{2}^{-1}, 0)^{2}\big]}{C_{\Psi}(H) ^{2}} -1 \right) \nonumber\\
&=& \frac{ M_{1}^{2H+1}M_{2}^{2H+1}}{C_{\Psi} (H)^{2}} \left[2^{2N (2H+1)} \mathbf{E}\Big[\widetilde{A_{M_{1}}(\ell_0, N) }^{2}\widetilde{A_{M_{2}}(\ell_0, N) }^{2} \Big] - \mathbf{E}\big[ c(M_{1}^{-1},0) ^{2}c(M_{2}^{-1}, 0)^{2}\big] \right. \nonumber \\
	&& - \left. \left( 2^{2N (2H+1)}\mathbf{E}\Big[\widetilde{A_{M_{1}}(\ell_0, N) }^{2}\Big]\mathbf{E}\Big[\widetilde{A_{M_{2}}(\ell_0, N) }^{2} \Big]- M_{1}^{-(2H+1)}M_{2}^{-(2H+1)} C_{\Psi} (H)^{2} \right) \right]. \nonumber \\ \label{4d-7}
\end{eqnarray}
First, we get, from Lemma \ref{l5},
\begin{eqnarray}
&&	\left|2^{2N (2H+1)} \mathbf{E}\Big[\widetilde{A_{M_{1}}(\ell_0, N) }^{2}\widetilde{A_{M_{2}}(\ell_0, N) }^{2} \Big] - \mathbf{E}\big[ c(M_{1}^{-1},0) ^{2}c(M_{2}^{-1}, 0)^{2}\big] \right| \leq C 2 ^{[ N ^{\beta}]\frac{ 2H-2}{q}}\nonumber \\ \label{4d-8}
\end{eqnarray}
while we have, using Lemma \ref{ll4},
\begin{eqnarray}
	&& \left| 2^{2N (2H+1)}\mathbf{E}\Big[\widetilde{A_{M_{1}}(\ell_0, N) }^{2}\Big]\mathbf{E}\Big[\widetilde{A_{M_{2}}(\ell_0, N) }^{2} \Big]- M_{1}^{-(2H+1)}M_{2}^{-(2H+1)} C_{\Psi} (H)^{2} \right| \nonumber\\
	&& \leq \left| 2^{N (2H+1)}\mathbf{E}\Big[\widetilde{A_{M_{1}}(\ell_0, N) }^{2}\Big]- M_{1}^{-(2H+1)}C_{\Psi} (H) \right|  \left| 2^{N (2H+1)}\mathbf{E}\Big[\widetilde{A_{M_{2}}(\ell_0, N) }^{2}\Big] \right| \nonumber \\
	&& + \left| M_{1}^{-(2H+1)}C_{\Psi} (H) \right| \left| 2^{N (2H+1)}\mathbf{E}\Big[\widetilde{A_{M_2}(\ell_0, N) }^{2}\Big]- M_{2}^{-(2H+1)}C_{\Psi} (H) \right| \nonumber \\
	&& \leq  C 2 ^{[ N ^{\beta}]\frac{ 2H-2}{q}} \label{4d-9}
\end{eqnarray}
We conclude the proof by plugging \eqref{4d-8} and \eqref{4d-9} into \eqref{4d-7} to obtain the inequality \eqref{8d-1}.  \qed  

For the proof of our main results, we will use the notion of { \it strong independence.} introduced in \cite{Tu}. 

\begin{definition}
Two square integrable random variables $F$ and $G$ admiting the chaos expansion $F=\sum _{n\geq 0} I_{n} (f_{n}), G= \sum _{n\geq 0} I_{n} (g_{n})$ where  $f_{n}, g_{n} \in L^2(\R ^{n})$ are symmetric for every $n\geq 0$, are strongly independent if for every $m,n\geq 0$, the random variables $ I_{n}(f_{n})$ and $ I_{m} (g_{m}) $ are independent. 
\end{definition}

We need two technical lemmas concerning the strong independent random variables in Wiener chaos. The first lemma concerns the strong independence of squares of independent chaotic random variables.

\begin{lemma} \cite[Lemma 2]{BoDiTu}\label{ll5}
	Consider two integers $p, q\geq 1$. Let $ F=I_{p} (f)$ and $G=I_{q} (g)$ with $ f\in L^2(\R ^{p}), g \in L^2(\R^{q})$. Assume that $F$ and $G$ are independent random variables. Then $F^2$ and $G^2$ are strongly independent.
\end{lemma}

The second lemma deals with scalar product of Malliavin derivative of strongly independent random variables. Basic notions of Malliavin calculus, together with the notation used in the following lemma are recalled in the Appendix.

\begin{lemma}\cite[Lemma 1]{Tu} and \cite[Lemma 4]{BoDiTu}\label{ll6}
	Assume that $F$ and $G$ are two strongly independent random variables in $\mathbb{D} ^{1,2}$. Then
	
	\begin{enumerate}
		\item 	\begin{equation*}
			\langle DF, D(-L) ^{-1}G\rangle =\langle D(-L) ^{-1} F, DG\rangle =0.
		\end{equation*}
		
		\item $\langle DF, D(-L) ^{-1}F\rangle $ and $\langle DG, D(-L) ^{-1}G\rangle$  are independent random variables. 
	\end{enumerate}
\end{lemma}

Let us now prove the CLT for the sequence (\ref{v1}) when  $M=1,...,d$ is fixed. In the sequel, $d_W$ stands for the Wasserstein distance. If $F$ and $G$ are two random variables, it is defined as
\[ d_W(F,G) = \sup_{h \in \mathcal{A}} | \mathbf{E}[h(F)]-\mathbf{E}[h(G)]|,\]
where $\mathcal{A}$ is the set of Lipschitz function $h:  \R \to \R$ such that
\[ \sup_{x,y \in \R, x \neq y} \frac{|h(x)-h(y)|}{|x-y|} \leq 1.\]
We refer to \cite[Appendix C1]{NP-book} for a comprehensive view on distances between probability distributions.

\begin{prop}\label{pp3}
Fix $ M=1,...,d$ and consider the sequence $ (V_{N, M, 1}, N \geq 1)$ defined by (\ref{v1}). Let us set	\begin{equation}
	\label{Mfh}
	K_{\Psi, H}= \frac{ \mathbf{E}  \big[c(1,0) ^{4}\big]}{ C_{\Psi}(H)^{2}}-1= \frac{ \mathbf{E}\big[ c(1,0) ^{4}\big]}{ (\mathbf{E}\big[ c(1,0) ^{2}\big]) ^{2}}-1.
\end{equation}
We have
\begin{equation*}
	V_{N,M,1} \underset{N \to \infty}{\longrightarrow^{(d)}} N (0, K_{\Psi, H}),
\end{equation*}
where $\to ^{(d)} $ stands for the convergence in distribution. Moreover, for $N$ large enough, 
\begin{equation*}
	d_{W} \left( V_{N, M, 1}; N (0, K_{\Psi, H})\right)\leq C2 ^{-\frac{[N ^{\gamma}]}{2}}.
\end{equation*}
\end{prop}
\noindent {\bf Proof: } From Theorem \ref{clt} and equation \eqref{cltbis} in the Appendix, we know that
\begin{eqnarray}
		d_{W} \left( V_{N, M, 1}; N (0, K_{\Psi, H})\right)&\leq& C \left( \left|\mathbf{E}  \big[V_{N, M,1} ^{2}\big]-	K_{\Psi, H} \right| + \sqrt{ \text{Var} \left[ \langle DV_{N, M, 1}, D(-L) ^{-1}V_{N, M, 1}\rangle\right] }\right)\nonumber \\
		&\leq & C \left( 2 ^{ [N ^{\beta}]\frac{2H-2}{q}}+ \sqrt{ \text{Var} \left( \langle DV_{N, M, 1}, D(-L) ^{-1}V_{N, M, 1}\rangle\right) }\right),\label{5d-1}
\end{eqnarray} 
where we used the estimate in Proposition \ref{pp2}. We start by computing the quantity $ \langle DV_{N, M, 1}, D(-L) ^{-1}V_{N, M, 1}\rangle$.  We have by (\ref{v1}), 

\begin{equation*}
	DV_{N, M,1}=\frac{ 2 ^{N (2H+1)}}{C_{\Psi}(H)\sqrt{ \vert \mathcal{L}_{N, \gamma}\vert }}\sum_{ \ell \in \mathcal{L}_{N, \gamma}}D\left( \widetilde{A_{M}(\ell, N) }^{2} -\mathbf{E}\Big[\widetilde{A_{M}(\ell, N) }^{2}\Big]\right)
\end{equation*}
and
\begin{equation*}
	D(-L) ^{-1}V_{N,M,1}=\frac{ 2 ^{N (2H+1)}}{C_{\Psi}(H)\sqrt{ \vert \mathcal{L}_{N, \gamma}\vert }}\sum_{ \ell \in \mathcal{L}_{N, \gamma}}D(-L) ^{-1}\left( \widetilde{A_{M}(\ell, N) }^{2}- \mathbf{E} \Big[\widetilde{A_{M}(\ell, N) }^{2}\Big] \right).
\end{equation*}
Thus, we express
\begin{eqnarray*}
	\langle DV_{N,M,1}, D(-L) ^{-1} V_{N,M,1} \rangle&=&\frac{ 2 ^{2N (2H+1)}}{C_{\Psi}(H)^{2}  \vert \mathcal{L}_{N, \gamma}\vert }\\
	&&\sum_{\ell, j\in \mathcal{L}_{N, \gamma}}\langle D\widetilde{A_{M}(\ell, N) }^{2}, D(-L) ^{-1}\left( \widetilde{A_{M}(j, N) }^{2}- \mathbf{E} \Big[\widetilde{A_{M}(j, N) }^{2}\Big] \right)\rangle.
\end{eqnarray*}
If $\ell\not= j$, then $\widetilde{A_{M}(\ell, N)}$ and $ \widehat{A_{M}(j, N)}$ are independent, by Lemma \ref{ll2}, point 1.  By Lemma \ref{ll5}, the random variables $\widetilde{A_{M}(\ell, N)}^{2} $ and $\widehat{A_{M}(j, N)}^{2}$ are strongly independent and the Lemma \ref{ll6} implies that 
\begin{equation*}
	\langle D\widehat{A_{M}(\ell, N) }^{2}, D(-L) ^{-1}\left( \widehat{A_{M}(j, N) }^{2}- \mathbf{E} \Big[\widehat{A_{M}(j, N) }^{2}\Big] \right)\rangle =0.
\end{equation*}
Therefore, we get
\begin{eqnarray*}
	&& \langle DV_{N,M,1}, D(-L) ^{-1} V_{N,M,1} \rangle \\&&=\frac{ 2 ^{2N (2H+1)}}{C_{\Psi}(H)^{2}  \vert \mathcal{L}_{N, \gamma}\vert } \sum_{\ell, \in \mathcal{L}_{N, \gamma}}\left\langle D\widehat{A_{M}(\ell, N) }^{2}, D(-L) ^{-1}\left( \widehat{A_{M}(\ell, N) }^{2}- \mathbf{E} \Big[\widehat{A_{M}(\ell, N) }^{2}\Big] \right)\right\rangle\\
	&&= \frac{ 2 ^{2N (2H+1)}}{C_{\Psi}(H)^{2}  \vert \mathcal{L}_{N, \gamma}\vert }
	\sum_{\ell, \in \mathcal{L}_{N, \gamma}} H_{\ell, N},
\end{eqnarray*}
where we used the notation
\begin{equation*}
	H_{\ell, N}= \left\langle D\widetilde{A_{M}(\ell, N) }^{2}, D(-L) ^{-1}\left( \widehat{A_{M}(\ell, N) }^{2}- \mathbf{E} \Big[\widehat{A_{M}(\ell, N) }^{2}\Big] \right)\right\rangle.
\end{equation*}
Next, we estimate the variance of $	\langle DV_{N, M,1}, D(-L) ^{-1} V_{N,M,1} \rangle$. For any $N\geq 1$, we can write
\begin{eqnarray*}
	&&	\text{Var} \left( 	\langle DV_{N,M,1}, D(-L) ^{-1} V_{N,M,1} \rangle\right)\\ 
	&=&\mathbf{E} \left[ 	\left(\langle DV_{N,M,1}, D(-L) ^{-1} V_{N,M,1} \rangle-\mathbf{E}\Big[	\langle DV_{N,M,1}, D(-L) ^{-1} V_{N,M,1} \rangle \Big]\right)^2\right]  \\
	&=&\frac{ 2 ^{4N (2H+1)}}{C_{\Psi}(H)^{4}  \vert \mathcal{L}_{N, \gamma}\vert ^{2} }\mathbf{E}\left[\sum_{\ell,j \in \mathcal{L}_{N, \gamma}} \left( H_{\ell, N}-\mathbf{E}\big[H_{\ell, N}\big]\right)  \left( H_{j, N}-\mathbf{E}\big[H_{j, N}\big]\right)\right].
\end{eqnarray*}
By the second part of Lemma \ref{ll6}, for every $\ell\not=j$ we have that $H_{\ell, N}$ and $ H_{j, N}$ are independent and so 
\begin{equation*}
	\mathbf{E}  \left[\left( H_{\ell, N}-\mathbf{E}\big[H_{\ell, N}\big]\right)  \left( H_{j, N}-\mathbf{E}\big[H_{j, N}\big]\right)\right] =0.
\end{equation*}
Consequently,
\begin{equation}\label{9n-4}
	\text{Var} \left[ 	\langle DV_{N,M,1}, D(-L) ^{-1} V_{N,M,1} \rangle\right] =\frac{ 2 ^{4N (2H+1)}}{C_{\Psi}(H)^{4}  \vert \mathcal{L}_{N, \gamma}\vert ^{2} }\sum_{\ell \in \mathcal{L}_{N, \gamma}} \mathbf{E} \left[\left( H_{\ell, N}-\mathbf{E}\big[H_{\ell, N}\big]\right)^{2} \right].
\end{equation}
We claim that for every $N\geq 1$ and $\ell \in \mathcal{L}_{N, \gamma}$ we have 
\begin{equation}\label{9n-1}
	\mathbf{E} \left[\left( H_{\ell, N}-\mathbf{E}\big[H_{\ell, N}\big]\right)^{2} \right]\leq C(q, H) 2 ^{-4N (2H+1)}<\infty,
\end{equation}
with $ C(q, H)>0$ not depending on $N, \ell$. To prove the inequality (\ref{9n-1}), we will use the chaos expansion of $ H_{\ell, N}$. Recall that $ \widetilde{A_{M}(\ell, N)}= I_{q} (\widehat{g}_{M}(\ell, N)$ with $ \widehat{g}_{M}(\ell, N)$ from (\ref{gnt}).  This implies, via the product formula (\ref{prod}),
\begin{equation*}
	\widetilde{A_{M}(\ell, N)}^{2}-\mathbf{E}\Big[\widetilde{A_{M}(\ell, N)}^{2}\Big] =\sum_{r=0} ^{q-1}r! \binom{q}{r}  ^{2}I_{2q-2r}\left( \widehat{g}_{M}(\ell, N)\otimes _{r} \widehat{g}_{M}(\ell, N)\right),
\end{equation*}
and
\begin{eqnarray*}
	&&	D_{\ast} 	\widetilde{A_{M}(\ell, N)}^{2}=\sum_{r=0} ^{q-1}r! \binom{q}{r} ^{2}(2q-2r)I_{2q-2r-1}\left( (\widehat{g}_{M}(\ell, N)\otimes _{r} \widehat{g}_{M}(\ell, N)(\cdot, \ast)\right),\\
	&& D_{\ast} (-L) ^{-1} \left( 	\widetilde{A_{M}(\ell, N)}^{2}-\mathbf{E}\Big[\widetilde{A_{M}(\ell, N)}^{2}\Big] \right) =\sum_{r=0} ^{q-1}r! \binom{q}{r} ^{2}I_{2q-2r-1}\left( (\widehat{g}_{M}(\ell, N)\otimes _{r} \widehat{g}_{M}(\ell, N)(\cdot, \ast)\right).
\end{eqnarray*}
As a consequence of the above two relations, we write

\begin{eqnarray}
	H_{\ell, N} &=& \sum_{ r_{1}, r_{2}=0} ^{q-1} r_{1} ! r_{2} !\binom{q}{r_1} ^{2}   \binom{q}{r_2}^{2} (2q-2r_{1})\nonumber \\
	&&\int_{\mathbb{R}} dx I_{2q-2r-1}\left( (\widehat{g}_{M}(\ell, N)\otimes _{r_{1}} \widehat{g}_{M}(\ell, N)(\cdot, x)\right)I_{2q-2r-1}\left( (\widehat{g}_{M}(\ell, N)\otimes _{r_{2}} \widehat{g}_{M}(\ell, N)(\cdot, x)\right)dx \nonumber\\
	&=&
	\sum_{ r_{1}, r_{2}=0} ^{q-1} r_{1} ! r_{2} ! \binom{q}{r_1} ^{2}  \binom{q}{r_2}^{2} (2q-2r_{1})\nonumber \\
	&& \sum_{a=0} ^{(2q-2r_{1}-1) \wedge (2q-2r_{2}-1)}a! \binom{2q-2r_{1}-1}{a} \binom{2q-2r_{2}-1}{a} \nonumber\\
	&&I _{4q-2r_{1}-2r_{2}-2-2a} \left( (\widehat{g}_{M}(\ell, N)\widetilde{\otimes} _{r_{1}} \widehat{g}_{M}(\ell, N) \otimes _{a+1} (\widehat{g}_{M}(\ell, N)\widetilde{\otimes} _{r_{2}} \widehat{g}_{M}(\ell, N)\right).\label{9n-2}
\end{eqnarray}
We remark, thanks to (\ref{4d-2}), that 
\begin{equation*}
	q! \Vert \widehat{g}_{M}(\ell, N)\Vert ^{2} _{ L ^{2}(\R ^{q})}= \mathbf{E} \Big[\widetilde{A_{M}(\ell, N)}^{2}\big] \leq C 2 ^{-N (2H+1)},
\end{equation*}
and, as a consequence, for every $r_{1}, r_{2}=0,..., q-1$ and for every $a=0,...,(2q-2r_{1}-1) \wedge (2q-2r_{2}-1)$, we have 
\begin{eqnarray}
	&&	\mathbf{E} \left[ I _{4q-2r_{1}-2r_{2}-2-2a} \left( (\widehat{g}_{M}(\ell, N)\widetilde{\otimes} _{r_{1}} \widehat{g}_{M}(\ell, N) )\otimes _{a+1} (\widehat{g}_{M}(\ell, N)\widetilde{\otimes} _{r_{2}} \widehat{g}_{M}(\ell, N))\right)^{2}\right] \nonumber\\
	&=&C (q, r_{1}, r_{2}, a) \Vert  (\widehat{g}_{M}(\ell, N)\widetilde{\otimes} _{r_{1}} \widehat{g}_{M}(\ell, N) \otimes _{a+1} (\widehat{g}_{M}(\ell, N)\widetilde{\otimes} _{r_{2}} \widehat{g}_{M}(\ell, N))\Vert ^{2} _{ M ^{2}(\mathbb{R}^{ 4q-2r_{1}-2r_{2}-2a-2})}\nonumber\\
	&\leq & C (q, r_{1}, r_{2}, a) \Vert  \widehat{g}_{M}(\ell, N)\widetilde{\otimes} _{r_{1}} \widehat{g}_{M}(\ell, N)\Vert ^{2} _{ L^2(\R ^{2q-2r_{1}})}  \Vert  \widehat{g}_{M}(\ell, N)\widetilde{\otimes} _{r_{2}} \widehat{g}_{M}(\ell, N)\Vert ^{2} _{ L^2(\R ^{2q-2r_{2}})}\nonumber\\
	&\leq &  C (q, r_{1}, r_{2}, a) \Vert \widehat{g}_{M}(\ell, N)\Vert ^{8} _{ L ^{2}(\R ^{q})}\nonumber\\
	&\leq &  C (q, r_{1}, r_{2}, a)2 ^{-4N (2H+1)}.\label{9n-3}
\end{eqnarray}

By combining the above estimates (\ref{9n-3}) and (\ref{9n-2}), we obtain the claim  (\ref{9n-1}). Coming back to (\ref{9n-4}), we get, using \eqref{eq:cardiup}, 
\begin{eqnarray}
	\text{Var} \left[ 	\langle DV_{N,M,1}, D(-L) ^{-1} V_{N, M,1} \rangle\right] &\leq& C(q,H, \Psi)\frac{ 2 ^{4N (2H+1)}}{  \vert \mathcal{L}_{N, \gamma}\vert ^{2} }\sum_{\ell \in \mathcal{L}_{N, \gamma}}2 ^{-4N (2H+1)}\nonumber\\
	&\leq& C(q,H, \Psi)2 ^{[N ^{\gamma}]}.\label{5d-2}
\end{eqnarray}
The conclusion is obtained via (\ref{5d-1}) and (\ref{5d-2}). \qed

We immediately deduce the asymptotic normality of the modified wavelet variation and its rate of convergence under the Wasserstein distance. 

\begin{theorem}\label{tt1}
	Consider the sequence $V_{N, M}$ defined by (\ref{vnl}). Then 
	\begin{equation*}
		V_{N,M} \underset{N \to \infty}{\longrightarrow^{(d)}} N (0, K_{\Psi, H}),
	\end{equation*}
	and for $N$ large enough, 
	\begin{equation*}
		d_{W} \left( V_{N, M}; N (0, K_{\Psi, H})\right)\leq C2 ^{-\frac{[N ^{\gamma}]}{2}}.
	\end{equation*}
\end{theorem}
	\noindent {\bf Proof: } The proof is a consequence of the results proven in Propositions \ref{pp1},  and \ref{pp3}. Indeed, by the triangle's inequality,
	\begin{eqnarray*}
			d_{W} \left( V_{N, M}; N (0, K_{\Psi, H})\right)&\leq &d_{W}  \left( V_{N, M,1}; N (0, K_{\Psi, H})\right)+ d_{W} (V_{N, M, 1}, V_{N, M})\\
			&\leq & d_{W}  \left( V_{N, M,1}; N (0, K_{\Psi, H})\right)+ \mathbf{E} \big[\vert V_{N, M, 1}- V_{N, M} \vert \big] \\
			&\leq & d_{W}  \left( V_{N, M,1}; N (0, K_{\Psi, H})\right)+ \mathbf{E} \big[\vert V_{N, M, 2}\vert \big] +\mathbf{E} \big[\vert V_{N, M, 3}\vert \big].  
	\end{eqnarray*}
	It suffices to use Proposition \ref{pp3} to bound the first summand in the right-hand side above and Proposition \ref{pp1} to bound the second and third summands. \qed

Let us now state and prove a multidimensional CLT for the modified wavelet variation. 

\begin{theorem}\label{thm:2}
	Let $V_{N, M}$ be given by (\ref{vnl}). Then the $d$-dimensional random vector $(V_{N, M}, M=1,,,,,d)$ converges in distribution, as $N \to\ \infty$, to the $d$-dimensional Gaussian vector $N(0, K)$, where $K=(K_{M_{1}, M_{2}})_{ M_{1}, M_{2}=1,...,d}$ is such that 
	\begin{equation}\label{cll1}
		K_{M, M}= K_{\Psi, H} \mbox{ for every } M=1,..., d
	\end{equation}
and 
\begin{equation}\label{cll2}
	K_{M _{1}, M_{2}}= \frac{(M_{1}M_{2})^{2H+1}\mathbf{E}\big[ c(M_{1}^{-1},0) ^{2}c(M_{2}^{-1}, 0)^{2}\big]}{C_{\Psi}(H) ^{2}} -1  \mbox{ for } M_{1}, M_{2}=1,...,d, M_{1}\not=M _{2}.
\end{equation}
\end{theorem}
\noindent {\bf Proof: } The triangle's inequality gives
\begin{eqnarray}
	d_{W} \left( ( V_{N, M}, M=1,...,d), N (0, K)\right) &\leq& 	d_{W} \left( ( V_{N, M,1}, M=1,...,d), N (0, K)\right)\nonumber \\
	&&+ \sum_{M=1} ^{d} \left( \mathbf{E} \big[\vert V_{N, M, 2}\vert\big] + \mathbf{E} \big[\vert V_{N, M, 3} \vert \big] \right).\label{8d-4}
\end{eqnarray}
By Theorem \ref{clt} and equation \eqref{cltbis} in the Appendix, we have 
\begin{eqnarray}
		d_{W} \left( ( V_{N, M,1}, M=1,...,d), N (0, K)\right) &\leq& C\left( \sqrt{ \sum_{M_{1}, M_{2}=1} ^{d} \left( \mathbf{E}\big[ V_{N, M_{1}, 1} V_{N, M_{2}, 1}\big]-K_{M_{1}, M_{2}} \right) ^{2} }\right.\nonumber\\
		&&\left.  + \sqrt{ \sum_{M_{1}, M_{2}=1} ^{d}\text{Var}\left[ \langle DV_{N, M_{1}, 1}, D(-L) ^{-1}V_{N, M_{2},1}\rangle \right] }\right) \nonumber \\\label{8d-3}
\end{eqnarray}

We need to compute the variance of the quantity $\langle D V_{N, M_{1}, 1}, D(-L) ^{-1} V_{N, M_{2}, 1} \rangle$. The situation when $ M_{1} = M_{2}$ has been treated in the proof of Proposition  \ref{pp3} and the case $M_{1} \not=M_{2}$ follows in a similar way. Actually, we express

\begin{eqnarray*}
	\langle DV_{N,M_{1},1}, D(-L) ^{-1} V_{N,M_{2},1} \rangle&=&\frac{ 2 ^{2N (2H+1)}}{C_{\Psi}(H)^{2}  \vert \mathcal{L}_{N, \gamma}\vert }\\
	&=& \frac{ 2 ^{2N (2H+1)}}{C_{\Psi}(H)^{2}  \vert \mathcal{L}_{N, \gamma}\vert }
	\sum_{\ell, \in \mathcal{L}_{N, \gamma}} H_{\ell, N}(M_{1}, M_{2})
\end{eqnarray*}
where we used the notation
\begin{equation*}
	H_{\ell, N}(M_{1}, M_{2})= \left\langle D\widehat{A_{M_{1}}(\ell, N) }^{2}, D(-L) ^{-1}\left( \widehat{A_{M_{2}}(\ell, N) }^{2}- \mathbf{E}\Big[ \widehat{A_{M_{2}}(\ell, N) }^{2}\Big] \right)\right\rangle.
\end{equation*}
This gives, as in (\ref{9n-4}),
\begin{eqnarray}\label{6d-5}
&&	\text{Var} \left[ 	\langle DV_{N,M_{1},1}, D(-L) ^{-1} V_{N,M_{2},1} \rangle\right] \\
&&=\frac{ 2 ^{4N (2H+1)}}{C_{\Psi}(H)^{4}  \vert \mathcal{L}_{N, \gamma}\vert ^{2} }\sum_{\ell \in \mathcal{L}_{N, \gamma}} \mathbf{E}\left[ \left( H_{\ell, N}(M_{1}, M_{2})-\mathbf{E}\big[H_{\ell, N}(M_{1}, M_{2})\big]\right)^{2}\right]\nonumber
\end{eqnarray}
and we can show that
\begin{equation*}
	\mathbf{E} \left[\left( H_{\ell, N}(M_{1}, M_{2}) -\mathbf{E}\big[H_{\ell, N}(M_{1}, M_{2})\big]\right)^{2}\right]\leq C(q, H) 2 ^{-4N (2H+1)}.
\end{equation*}
We obtain, from \eqref{eq:cardiup},
\begin{equation}\label{8d-2}
		\text{Var} \left[ 	\langle DV_{N,M_{1},1}, D(-L) ^{-1} V_{N,M_{2},1} \rangle\right]\leq C 2 ^{-[ N ^{\gamma}]}.
\end{equation}
By combining (\ref{8d-1}), (\ref{8d-2})  and (\ref{8d-3}), we conclude
\begin{equation*}
		d_{W} \left( ( V_{N, M,1}, M=1,...,d), N (0, K)\right) \leq C 2 ^{ -\frac{ [ N^{\gamma}]}{2}}.
\end{equation*}
We finish the proof by plugging this estimate, together with the  estimates in Proposition \ref{pp1}, into (\ref{8d-4}).  \qed 

\section{Discretization of the wavelet variation}\label{sec:discret}
We will define an estimator for the Hurst index of the Hermite process based on the  modified wavelet variation $V_{N, M}$ given by (\ref{vnl}). In order to obtain an expression of the estimator which can be numerically computed, we need to discretize the wavelet coefficient $ A_{M}(\ell, N)$.

We will consider the following Riemann approximation: for $M=1,..., d$, $N\geq 1$ and for $\ell \in \mathcal{L}_{N, \gamma}$, we set 
\begin{equation}
	\label{eln}
E_{M} (\ell, N)= \sqrt{ \frac{1}{ M 2 ^{N}}}\frac{1}{ 2 ^{N}} \sum _{k=1} ^{2 ^{N}}\Psi \left( \frac{k}{2 ^{N}}\right) Z _{ \frac{ k2 ^{-N} + M\ell 2 ^{[N ^{\beta}]}}{ M 2 ^{N}}},
\end{equation}
and
\begin{equation}
	\label{hatv}
	\widehat{V}_{N. M}= \frac{1}{ \sqrt{ \vert \mathcal{L}_{ N, \gamma} \vert }}\sum _{\ell \in \mathcal{L}_{ N, \gamma}}\left[ \frac{ E_{M} (\ell, N) ^{2}}{ \mathbf{E} \big[A_{M} (\ell, N) ^{2}\big]}-1 \right]. 
\end{equation}
We will prove that $ V_{N, M}- \widehat{V}_{N, M}$ converges to zero in $L ^{1} (\Omega)$ as $ N \to \infty$ and we will deduce that the discretized wavelet variation $ \widehat{V}_{N, M}$ also satisfies a CLT. We start by evaluating the difference between $ A_{M} (\ell, N)$ and its discretized counterpart $ E_{M} (\ell, N)$.

\begin{prop}\label{pp4}
Let us denote
\begin{equation}
	\label{tnl}
	t_{N,M}=\frac{ A_{M}(\ell, N)}{ \sqrt{ \mathbf{E}\big[A_{M} (\ell, N) ^{2}\big]}}- \frac{ E_{M}(\ell, N)}{ \sqrt{ \mathbf{E} \big[A_{M} (\ell, N) ^{2}\big]}}, \hskip0.4cm N\geq 1, M=1,...,d.
\end{equation}
Then, for $N$ large enough, we have
\begin{equation*}
	\mathbf{E} \big[t_{N, M} ^{2}\big] \leq C 2 ^{-N}2 ^{ [ N ^{\beta}]2H}.
\end{equation*}
\end{prop}
\noindent {\bf Proof: } By (\ref{4d-3}), we get 
\begin{equation*}
	t_{N, M}= \sqrt{ \frac{ M ^{2H+1}}{C _{\Psi }(H)} 2 ^{N (2H+1)}}( A_{M} (\ell, N)- E _{M} (\ell, N)),
\end{equation*}
and then 
\begin{equation*}
	\mathbf{E}\big[ t_{N, M} ^{2}\big]= \frac{ M ^{2H+1}}{C _{\Psi }(H)} 2 ^{N (2H+1)}( \mathbf{E} \big[A_{M} (\ell, N) ^{2}\big]- 2\mathbf{E} \big[A_{M} (\ell, N) E_{M} (\ell, N)\big] + \mathbf{E} \big[E_{M} (\ell, N) ^{2}\big]).
\end{equation*}
By (\ref{al}) and (\ref{eln}), we write
\begin{eqnarray*}
&&	\mathbf{E} \big[A_{M} (\ell, N) E_{M} (\ell, N)\big]\\
&=&\frac{1}{ M 2 ^{2N}}\sum_{k=1} ^{2 ^{N}}\Psi \left( \frac{k}{ 2 ^{N}}\right) \int_{\mathbb{R}} dx \Psi (x) \mathbf{E} \left[ Z _{ \frac{ x+ M\ell 2 ^{[N ^{\beta}]}}{ M 2 ^{N}}}Z _{ \frac{ k2 ^{-N} + M\ell 2 ^{[N ^{\beta}]}}{ M 2 ^{N}}}\right] \\
	&=& \frac{1}{ M 2 ^{2N}}\frac{1}{2 (M 2 ^{N})^{2H}}\sum_{k=1} ^{2 ^{N}}\Psi \left( \frac{k}{ 2 ^{N}}\right) \int_{\mathbb{R}} dx \Psi (x)\left[ (x+M\ell 2 ^{[N ^{\beta}]})^{2H} -\vert x-k2 ^{-N}\vert ^{2H}\right],
\end{eqnarray*}
where we used the assumption (\ref{mom1}). Also, we have
\begin{eqnarray*}
	\mathbf{E} \big[E_{M} (\ell, N) ^{2}\big] &=& \frac{1}{ M 2 ^{3N}}\sum_{k,j=1} ^{2 ^{N}}\Psi \left( \frac{k}{ 2 ^{N}}\right) \Psi \left( \frac{j}{ 2 ^{N}}\right) \mathbf{E} \left[ Z _{ \frac{ k 2 ^{-N}+ M\ell 2 ^{[N ^{\beta}]}}{ M 2 ^{N}}} Z _{ \frac{ j 2 ^{-N}+ M\ell 2 ^{[N ^{\beta}]}}{ M 2 ^{N}}}\right]\\
	&=& \frac{1}{ M 2 ^{3N}}\frac{1}{2 (M 2 ^{N})^{2H}}\sum_{k,j=1} ^{2 ^{N}}\Psi \left( \frac{k}{ 2 ^{N}}\right) \Psi \left( \frac{j}{ 2 ^{N}}\right) \\
	&&\left[  (k 2 ^{-N}+M\ell 2 ^{[N ^{\beta}]})^{2H}+(j 2 ^{-N}+M\ell 2 ^{[N ^{\beta}]})^{2H}-\vert k2 ^{-N}-j2 ^{-N} \vert ^{2H} \right].
\end{eqnarray*}
Consequently, we express
\begin{eqnarray*}
	\mathbf{E} \big[t_{N, M} ^{2}\big]&=& \frac{ M ^{2H+1}}{ C_{\Psi}(H)} \left[ M ^{-(2H+1)}C_{\Psi }(H) - \right.\\
	&&\left.M ^{-2H-1} \frac{1}{ 2 ^{N}} \sum_{k=1} ^{2 ^{N}}\Psi \left( \frac{k}{ 2 ^{N}}\right) \int_{\mathbb{R}} dx \Psi (x)\left[ (x+M\ell 2 ^{[N ^{\beta}]})^{2H} -\vert x-k2 ^{-N}\vert ^{2H}\right]\right.\\
	&&\left.+ \frac{1}{2} M ^{-2H-1} \frac{1}{ 2 ^{2N}}\sum_{k,j=1} ^{2 ^{N}}\Psi \left( \frac{k}{ 2 ^{N}}\right) \Psi \left( \frac{j}{ 2 ^{N}}\right)\right. \\
	&&\left. \left[  (k 2 ^{-N}+M\ell 2 ^{[N ^{\beta}]})^{2H}+(j 2 ^{-N}+M\ell 2 ^{[N ^{\beta}]})^{2H}-\vert k2 ^{-N}-j2 ^{-N} \vert ^{2H} \right] \right]\\
	&=&1+ \frac{ 1}{ C_{\Psi}(H)}\frac{1}{2^{N}}\sum_{k=1} ^{2 ^{N}}\Psi \left( \frac{k}{ 2 ^{N}}\right) \int_{\mathbb{R}} dx \Psi (x)\vert x-k 2 ^{-N}\vert ^{2H}\\
	&&- \frac{1}{2}\frac{ 1}{ C_{\Psi}(H)}\frac{1}{2^{2N}}\sum_{k,j=1} ^{2 ^{N}}\Psi \left( \frac{k}{ 2 ^{N}}\right) \Psi \left( \frac{j}{ 2 ^{N}}\right)\vert k2 ^{-N}-j2 ^{-N}\vert ^{2H}\\
	&&-C_{\Psi }(H) ^{-1} \left[ \frac{1}{ 2 ^{N}} \sum_{k=1} ^{2 ^{N}}\Psi \left( \frac{k}{ 2 ^{N}}\right) \int_{\mathbb{R}} dx \Psi (x)(x+M\ell 2 ^{[N ^{\beta}]})^{2H} \right. \\
	&&\left. - \frac{1}{ 2 ^{2N}}\sum_{k,j=1} ^{2 ^{N}}\Psi \left( \frac{k}{ 2 ^{N}}\right) \Psi \left( \frac{j}{ 2 ^{N}}\right) (k 2 ^{-N}+M\ell 2 ^{[N ^{\beta}]})^{2H}\right] \\
	&=:& 1+ C_{\Psi } (H) ^{-1} (T_{1,N}-\frac{1}{2} T_{2, N} )-C_{\Psi }(H) ^{-1} Q_{N},
\end{eqnarray*}
with
\begin{eqnarray*}
	T_{1,N}=\frac{1}{2 ^{N}}\sum_{k=1} ^{2 ^{N}}\Psi \left( \frac{k}{ 2 ^{N}}\right) \int_{\mathbb{R}} dx \Psi (x)\vert x-k 2 ^{-N}\vert ^{2H},
\end{eqnarray*}

\begin{equation*}
	T_{2, N}=\frac{1}{2^{2N}}\sum_{k,j=1} ^{2 ^{N}}\Psi \left( \frac{k}{ 2 ^{N}}\right) \Psi \left( \frac{j}{ 2 ^{N}}\right)\vert k2 ^{-N}-j2 ^{-N}\vert ^{2H},
\end{equation*}
and
\begin{eqnarray*}
	Q_{N} &=&  2 ^{-N} \sum_{k=1} ^{2 ^{N}} \Psi \left( \frac{k}{2 ^{N}}\right) \left[ \int_{\mathbb{R}} dx \Psi (x) (x+\ell M 2 ^{ [n ^{\beta}]})^{2H} -\right.\\
	&&\left.2 ^{-N} \sum_{j=1} ^{2 ^{N}}  \Psi \left( \frac{j}{2 ^{N}}\right) (k 2 ^{-N}+\ell M 2 ^{ [n ^{\beta}]})^{2H}\right]
\end{eqnarray*}

We show that 
\begin{equation*}
	T_{1,N} \underset{N \to \infty}{\longrightarrow} \int_{0} ^{1} \int _{0} ^{1}  \Psi (x) \Psi (y) \vert x-y\vert ^{2H-2}dxdy = -2C_{\Psi}(H)
\end{equation*}
and
\begin{equation}\label{6d-1}
	\vert T_{1, N}- (-2C_{\Psi}(H))\vert \leq C 2 ^{-N}.
\end{equation}
We can decompose $ T_{1,N}$ as follows. 
\begin{eqnarray*}
	T_{1,N}&=&  \int_{0} ^{1} \int _{0} ^{1}  \Psi (x) \Psi (y) \vert x-y\vert ^{2H-2}dxdy\\
	&&+ \sum_{k=1} ^{2 ^{N}} \int_{ \frac{ k-1} {2 ^{N}}}^{\frac{ k}{2 ^{N}}}dy \left( \Psi \left( \frac{k}{ 2 ^{N}}\right)- \Psi (y) \right) \int_{\mathbb{R}} dx \Psi (x) \vert x-k2 ^{-N}\vert ^{2H}\\
	&&+ \sum_{k=1} ^{2 ^{N}} \int_{ \frac{ k-1} {2 ^{N}}}^{\frac{ k}{2 ^{N}}}dy\Psi (y) \int_{\mathbb{R}} \Psi (x) dx \left( \vert x+ k 2 ^{-N}\vert ^{2H}- \vert x-y\vert ^{2H} \right)\\
	&=:& -2 C_{\Psi }(H) + R_{1, N}+ R _{2, N}.
\end{eqnarray*}
Let us bound the two rest terms, denoted by $ R_{1, N}$ and $R_{2,N}$. We have, since $ \Psi $ is of class $ C ^{1}$ with support contained in $[0,1]$, 
\begin{eqnarray}
	\vert R_{1, N}\vert &\leq & C \sum_{k=1} ^{2 ^{N}} \int_{ \frac{ k-1} {2 ^{N}}}^{\frac{ k}{2 ^{N}}}dy\vert y- k 2 ^{-N}\vert  \int_{\mathbb{R}} dx \vert \Psi (x)\vert  \vert x-k2 ^{-N}\vert ^{2H}\nonumber\\
	&\leq & C 2 ^{-N} \sum_{k=1} ^{2 ^{N}} \int_{ \frac{ k-1} {2 ^{N}}}^{\frac{ k}{2 ^{N}}}dy \int_{\mathbb{R}} dx \vert \Psi (x)\vert \leq C 2 ^{-N}.\label{6d-2}
\end{eqnarray}
Concerning the summand $R_{2,N}$, let us remark that the function \[f:[0,1]\to \mathbb{R} \, : \, a \mapsto \vert x-a\vert ^{2H}\] satisfies $ \vert f'(a)\vert \leq C $ as soon as $H>\frac{1}{2}$. Thus, we get 
\begin{eqnarray}
	\vert R_{2, N} \vert &\leq & C  \sum_{k=1} ^{2 ^{N}} \int_{ \frac{ k-1} {2 ^{N}}}^{\frac{ k}{2 ^{N}}}dy\vert \Psi (y)\vert  \vert y-k 2 ^{-N}\vert \int_{\mathbb{R}} \vert \Psi (x) \vert dx\nonumber \\
	&\leq & C 2 ^{-N} \sum_{k=1} ^{2 ^{N}} \int_{ \frac{ k-1} {2 ^{N}}}^{\frac{ k}{2 ^{N}}}dy\vert \Psi (y)\vert  \int_{\mathbb{R}} \vert \Psi (x) \vert dx \leq C 2^{-N}.\label{6d-3}
\end{eqnarray}
The estimates (\ref{6d-2}) and (\ref{6d-3}) imply (\ref{6d-1}).

In the same way, we can prove that $ T_{2, N} \underset{ N \to \infty}{\longrightarrow} -2C_{\Psi}(H)$ and 
\begin{equation*}
	\left| T_{2, N}- (-2C_{\Psi}(H) )\right| \leq C 2 ^{-N}.
\end{equation*}
We now show that $Q_{N}$ converges to zero as $N$ tends to infinity. One has 
\begin{eqnarray*}
	Q_{N} &=&  2 ^{-N} \sum_{k=1} ^{2 ^{N}} \Psi \left( \frac{k}{2 ^{N}}\right) \left[ \int_{\mathbb{R}} dx \Psi (x) (x+\ell M 2 ^{ [n ^{\beta}]})^{2H} -2 ^{-N} \sum_{j=1} ^{2 ^{N}}  \Psi \left( \frac{j}{2 ^{N}}\right) (j 2 ^{-N}+\ell M 2 ^{ [n ^{\beta}]})^{2H}\right]\\
	&=& 2 ^{-N}  \sum_{k=1} ^{2 ^{N}} \Psi \left( \frac{k}{2 ^{N}}\right)  \sum_{j=1} ^{2 ^{N}} \int_{ \frac{j-1} {2 ^{N}}}^{\frac{j}{ 2 ^{N}} }dx  \left[ \Psi (x) (x+\ell M 2 ^{ [n ^{\beta}]})^{2H}-\Psi \left( \frac{j}{2 ^{N}}\right) (j 2 ^{-N}+\ell M 2 ^{ [n ^{\beta}]})^{2H}\right].
\end{eqnarray*}
Hence
\begin{eqnarray}
	\vert Q_{N} \vert &\leq & 2 ^{-N}  \sum_{k=1} ^{2 ^{N}} \left| \Psi \left( \frac{k}{2 ^{N}}\right) \right| \sum_{j=1} ^{2 ^{N}} \int_{ \frac{j-1} {2 ^{N}}}^{\frac{j}{ 2 ^{N}} }dx \vert \Psi (x) \vert \left|  (x+\ell M 2 ^{ [n ^{\beta}]})^{2H}- (j 2 ^{-N}+\ell M 2 ^{ [n ^{\beta}]})^{2H} \right| \nonumber\\
	&&+2 ^{-N}  \sum_{k=1} ^{2 ^{N}} \left| \Psi \left( \frac{k}{2 ^{N}}\right) \right| \sum_{j=1} ^{2 ^{N}} \int_{ \frac{j-1} {2 ^{N}}}^{\frac{j}{ 2 ^{N}} }dx (j 2 ^{-N}+\ell M 2 ^{ [n ^{\beta}]})^{2H} \left| \Psi (x)-  \Psi \left( \frac{j}{2 ^{N}}\right)\right|.\nonumber
\end{eqnarray}
We use the bounds, for $ x\in (\frac{j-1}{ 2 ^{N}}, \frac{j}{2 ^{N}})$,
\begin{equation*}
	\vert (x+\ell M 2 ^{ [n ^{\beta}]})^{2H}- (j 2 ^{-N}+\ell M 2 ^{ [n ^{\beta}]})^{2H}\vert \leq C 2 ^{-N}2 ^{(2H-1) [N ^{\beta}]}
\end{equation*}
and
\begin{equation*}
\left| \Psi (x)-  \Psi \left( \frac{j}{2 ^{N}}\right)\right| \leq C 2 ^{-N}.
\end{equation*}
Thus, we reach
\begin{equation}
	\label{6d-6}
	\vert Q_{N}\vert \leq  C 2 ^{-N}2 ^{ [ N ^{\beta}]2H}.
\end{equation}
From (\ref{6d-1}) and (\ref{6d-6}), we obtained the conclusion. \qed

Next, let us compare the wavelet variation (\ref{vnl})  and its discretized version (\ref{hatv}). 

\begin{prop}\label{pp5}
	For $N$ sufficiently large, we have
	\begin{equation*}
		\mathbf{E} \big[\vert V_{N,M}- \widehat{V}_{N, M}\vert \big]\leq C 2 ^{-N}2 ^{ [ N ^{\beta}]2H} 2 ^{-\frac{ [N ^{\gamma}]}{2}}.
	\end{equation*}
\end{prop}
\noindent {\bf Proof: } By (\ref{vnl}) and (\ref{hatv}), we can write
\begin{eqnarray*}
	&&V_{N,M}-\widehat{V}_{N. M} = \frac{1}{ \sqrt{ \vert \mathcal{L}_{N, \gamma}\vert }}\sum_{\ell\in \mathcal{L}_{N, \gamma}}\\
	&&\left( \frac{A_{M}(\ell, N)}{ \sqrt{ \mathbf{E} \big[A_{M} (\ell, N) ^{2}\big]}}-\frac{E_{M}(\ell, N)}{ \sqrt{ \mathbf{E} \big[A_{M} (\ell, N) ^{2}\big]}}\right)\left( \frac{A_{M}(\ell, N)}{ \sqrt{ \mathbf{E} \big[A_{M} (\ell, N) ^{2}\big]}}+\frac{E_{M}(\ell, N)}{ \sqrt{ \mathbf{E} \big[A_{M} (\ell, N) ^{2}\big]}}\right)
\end{eqnarray*}
and thus, by the Cauchy-Schwarz's inequality,
\begin{eqnarray*}
	&&\mathbf{E} \left[\left| V_{N,M}-\widehat{V}_{N. M}\right|\right] \\
	&\leq & \frac{1}{ \sqrt{ \vert \mathcal{L}_{N, \gamma}\vert }}\sum_{\ell\in \mathcal{L}_{N, \gamma}}\\
	&&\mathbf{E} \left[ \left| \frac{A_{M}(\ell, N)}{ \sqrt{ \mathbf{E} A_{M} (\ell, N) ^{2}}}-\frac{E_{M}(\ell, N)}{ \sqrt{ \mathbf{E} A_{M} (\ell, N) ^{2}}}\right| \cdot  \left| \frac{A_{M}(\ell, N)}{ \sqrt{ \mathbf{E} A_{M} (\ell, N) ^{2}}}+\frac{E_{M}(\ell, N)}{ \sqrt{ \mathbf{E} A_{M} (\ell, N) ^{2}}}\right| \right]\\
	&\leq &  \frac{1}{ \sqrt{ \vert \mathcal{L}_{N, \gamma}\vert }}\sum_{\ell\in \mathcal{L}_{N, \gamma}}\\
	&&\left( \mathbf{E} \left[\left| \frac{A_{M}(\ell, N)}{ \sqrt{ \mathbf{E} A_{M} (\ell, N) ^{2}}}-\frac{E_{M}(\ell, N)}{ \sqrt{ \mathbf{E} A_{M} (\ell, N) ^{2}}}\right| ^{2} \right] \right) ^{\frac{1}{2}}\left( \mathbf{E} \left[\left| \frac{A_{M}(\ell, N)}{ \sqrt{ \mathbf{E} A_{M} (\ell, N) ^{2}}}+\frac{E_{M}(\ell, N)}{ \sqrt{ \mathbf{E} A_{M} (\ell, N) ^{2}}}\right| ^{2}\right] \right) ^{\frac{1}{2}}.
\end{eqnarray*}
Since, clearly, 
\begin{equation*}
	 \mathbf{E} \left[\left| \frac{A_{M}(\ell, N)}{ \sqrt{ \mathbf{E} A_{M} (\ell, N) ^{2}}}+\frac{E_{M}(\ell, N)}{ \sqrt{ \mathbf{E} A_{M} (\ell, N) ^{2}}}\right| ^{2} \right]\leq C,
\end{equation*}
we get from Proposition \ref{pp4},
\begin{eqnarray*}
	\mathbf{E}\left[ \left| V_{N,M}-\widehat{V}_{N. M}\right|\right] &\leq & \frac{1}{ \sqrt{ \vert \mathcal{L}_{N, \gamma}\vert }}\sum_{\ell\in \mathcal{L}_{N, \gamma}}2 ^{ [N ^{\beta}]H}2 ^{-\frac{N}{2}}\\
	&\leq &  C 2 ^{-\frac{N}{2}}2 ^{ [ N ^{\beta}]H} 2 ^{\frac{ [N ^{\gamma}]}{2}}\underset{N\to \infty}{\longrightarrow}0.
\end{eqnarray*}

The above result implies that the discretized wavelet variation has the same limit behavior in distribution as $ V_{N, M}$. 

\begin{corollary}\label{cor1}
	Let $\widehat{V}_{N, M}$ be given by (\ref{hatv}) Then the $d$-dimensional random vector $(\widehat{V}_{N, M}, M=1,,,,,d)$ converges in distribution, as $N \to\ \infty$, to the $d$-dimensional Gaussian vector $N(0, K)$, where the matrix $K$ is given by (\ref{cll1}) and (\ref{cll2}).
\end{corollary}
\noindent {\bf Proof: } The result follows from Theorem \ref{tt1} and Proposition \ref{pp5}. 

We finish this section with a short result concerning the almost-sure convergence to $0$ of the sequence $(\widehat{V}_{N, M}/\sqrt{ \vert \mathcal{L}_{N, \gamma} \vert }, N \geq 1)$. This fact will be particularly useful in the next section.

\begin{prop}\label{pp6}
	We have, for all $M \in \{1,\dots,d\}$,\begin{equation}
		\frac{\widehat{V}_{N, M}}{ \sqrt{ \vert \mathcal{L}_{N, \gamma} \vert }} \underset{N \to \infty}{\longrightarrow} 0 
	\end{equation}
almost surely, at the fast rate $2^{-N^a}$, where $a \in (0,\gamma)$ is arbitrary and fixed
\end{prop}
\noindent {\bf Proof: } 
Let us remark that Proposition \ref{pp5} and Theorem \ref{tt1} entail that, for $M=1,\dots,d$, the sequence $(\mathbf{E}\big[\vert \widehat{V}_{N, M}\vert \big])_N$ is bounded. Therefore, as a consequence of \eqref{eq:cardi}, we get, for such $M$,
\begin{equation}\label{eq:avantmarkov}
 \mathbf{E}\big[\vert \widehat{V}_{N, M}\vert \big] \leq C 2 ^{-\frac{[N ^{\gamma}]}{2}}
\end{equation}
The proofs is then obtained by a simple Borel-Cantelli argument. Let $a\in (0, \gamma)$ be arbitrary and fixed. We deduce from \eqref{eq:avantmarkov} and Markov inequality that, for $M=1,\dots,d$,
\begin{equation*}
\sum_{N\geq 1}	\mathbf{P} \left( \left| 	\frac{\widehat{V}_{N, M}}{ \sqrt{ \vert \mathcal{L}_{N, \gamma} \vert }}\right| \geq 2 ^{ N^{a}}\right)\leq \sum_{N\geq 1}  2 ^{ N^{a}} 2 ^{-\frac{[N ^{\gamma}]}{2}}\mathbf{E}\big[\vert \widehat{V}_{N, M}\vert \big]\leq C\sum_{N\geq 1}  2 ^{ N^{a}} 2 ^{-\frac{[N ^{\gamma}]}{2}}< \infty.
\end{equation*}
\qed

\section{Estimation of the Hurst parameter}\label{sec:estima}
Let us introduce the sequences

\begin{equation}
	\label{snl}
S_{N, M}= \frac{1}{ \vert \mathcal{L}_{N, \gamma}\vert } \sum_{\ell \in \mathcal{L}_{N, \gamma}}A_{M} (\ell, N) ^{2}
\end{equation}
and
\begin{equation}
	\label{hats}
	\widehat{S}_{N, M}= \frac{1}{ \vert \mathcal{L}_{N, \gamma}\vert } \sum_{\ell \in \mathcal{L}_{N, \gamma}}E_{M} (\ell, N) ^{2}
\end{equation}
where $A_{M}(\ell, N),  E_{M}(\ell, N)$ are the wavelet coefficients given by (\ref{al}) and  (\ref{eln}), respectively. We clearly have, via Proposition \ref{pp4}, 
\begin{equation*}
	\mathbf{E} \big[\vert A_{M}(\ell, N) - E_{M}  (\ell, N)\vert ^{2}\big] \leq C 2 ^{[N ^{\beta}]2H}2 ^{-N}\mathbf{E} \big[A_{M}(\ell, N) ^{2}\big] \leq C 2 ^{[N ^{\beta}]2H}2 ^{-N}2 ^{-N (2H+1)}
\end{equation*}
so, by Cauchy-Schwarz's inequality, we get
\begin{eqnarray*}
	\mathbf{E}\big[ \vert S_{N, M}- \widehat{S}_{N, M}\vert\big] &\leq &\frac{1}{ \vert \mathcal{L}_{N, \gamma}\vert } \sum_{\ell \in \mathcal{L}_{N, \gamma}}\mathbf{E} \big[\vert A_{M}(\ell, N)^2 - E_{M}  (\ell, N)^2\vert \big] \\
	&\leq &\frac{1}{ \vert \mathcal{L}_{N, \gamma}\vert } \sum_{\ell \in \mathcal{L}_{N, \gamma}}\left( 
	\mathbf{E}\big[ \vert A_{M}(\ell, N) - E_{M}  (\ell, N)\vert ^{2}\big] \right) ^{\frac{1}{2}}\left( \mathbf{E} \big[\vert A_{M}(\ell, N) + E_{M}  (\ell, N)\vert ^{2}\big] \right) ^{\frac{1}{2}}\\
	&\leq & C\sqrt{ 2 ^{[N ^{\beta}]2H}2 ^{-N}2 ^{-N (2H+1)}} \sqrt{ 2 ^{-N (2H+1)}},
\end{eqnarray*}
which entails
\begin{equation}\label{7d-1}
 S_{N, M} -\widehat{S}_{N, M} \underset{N \to \infty}{\longrightarrow} 0  \mbox{ in } L ^{1}(\Omega).  
\end{equation}

We also have, due to (\ref{4d-3}),
\begin{equation*}
	\mathbf{E} \big[S_{N, M}\big] = (M 2 ^{N}) ^{-(2H+1)}C_{\Psi }(H),
\end{equation*}
so
\begin{equation}\label{7d-2}
	\log \mathbf{E} \big[S_{N, M}\big]= -(2H+1) \log (M 2 ^{N})+ \log C_{\Psi }(H).
\end{equation}
To construct an estimator for the Hurst index of the Hermite process, we use the following standard procedure: in (\ref{7d-2}) we approximate $\mathbf{E} [S_{N, M}]$ by $ S_{N, M}$ and thus, due to (\ref{7d-1}), by $\widehat{S}_{N, M}$. We then have
\begin{equation*}
	\log \widehat{S}_{N, M}\sim -(2H+1) \log (M 2 ^{N})+ \log C_{\Psi }(H).
\end{equation*}
Next, we make a log-regression of $\left( \log \widehat{S}_{N, M}, M=1,...,d\right)$ on $\left( \log (M 2 ^{N}), M=1,...,d\right)$ which leads to the estimator $\widehat{H}_{N}$ given by
\begin{equation}
\widehat{H}_{N} =-\frac{1}{2}\frac{\left( \sum_{M=1}^{d} \log  \widehat{S}_{N,M} \log M  \right)-\left(\sum_{M=1}^{d} \log  \widehat{S}_{N,M} \right)\left(\sum_{M=1}^{d} \log M  \right) }{ \left(\sum_{M=1}^{d} (\log M )^2 \right)-\left(\sum_{M=1}^{d} \log M  \right)^2}-\frac{1}{2} \label{est1}
\end{equation}
We state and prove the limit behavior of the above estimator.

\begin{theorem}
	Let $ \widehat{H}_{N}$ be given by (\ref{est1}). 	Then the estimator $\widehat{H}_{N}$ is strongly consistent, i.e.
	\begin{equation}\label{almostsureconv}
		\widehat{H}_{N} \underset{N \to \infty}{\longrightarrow}  H \mbox{ almost surely}. 
	\end{equation}
	Moreover, 
	
	\begin{equation}\label{9d-2}
		 \sqrt{ \vert \mathcal{L}_{N, \gamma} \vert }( H-\widehat{H}_{N})\underset{N \to \infty}{\longrightarrow^{(d)}} N(0, \sigma ^{2}),
	\end{equation}
	with $\sigma ^{2} =\frac{1}{4}(L_d^T L_d)^{-1} L_d K L_d^T (L_d^T L_d)^{-1}$, where $K$ is the matrix defined in Theorem \ref{thm:2} and $L_d$ is the matrix with $(L_d)_{M,1}=\log M$ and $(L_d)_{M,2}=1$, for $M=1,\dots,d$.
\end{theorem}
\noindent {\bf Proof: } To obtain the asymptotic properties of the estimator (\ref{est1}), we use the limit theorems obtained for the discretized wavelet variation $\widehat{V}_{N, M}$. We notice the following link between $ \widehat{S}_{N, M}$ and $\widehat{V}_{N, M}$:

\begin{equation*}
	\frac{\widehat{V}_{N, M}}{ \sqrt{ \vert \mathcal{L}_{N, \gamma} \vert }}+1 = \frac{ (M2 ^{N}) ^{(2H+1)}}{C_{\Psi }(H)} \widehat{S}_{N, M},
\end{equation*}
and thus
\begin{equation*}
\log \widehat{S}_{N, M}= \log\left(\frac{\widehat{V}_{N, M}}{ \sqrt{ \vert \mathcal{L}_{N, \gamma} \vert }}+1 \right)-(2H+1) (\log M +N \log(2))-\log C_\Psi(H).
\end{equation*}
From this last equality, we deduce
\begin{eqnarray}
&&\widehat{H}_{N}-H = - \frac{1}{2}\frac{\left( \sum_{M=1}^{d} \log\left(\frac{\widehat{V}_{N, M}}{ \sqrt{ \vert \mathcal{L}_{N, \gamma} \vert }}+1 \right) \log M  \right)-\left(\sum_{M=1}^{d} \log\left(\frac{\widehat{V}_{N, M}}{ \sqrt{ \vert \mathcal{L}_{N, \gamma} \vert }}+1 \right) \right)\left(\sum_{M=1}^{d} \log M  \right) }{ \left(\sum_{M=1}^{d} (\log M )^2 \right)-\left(\sum_{M=1}^{d} \log M  \right)^2}  \nonumber \\ \label{9d-3}
\end{eqnarray}
The strong consistence \eqref{almostsureconv} is then a straightforward consequence of Proposition \ref{pp6}.

To prove \eqref{9d-2}, let us recall that the inequality 
\begin{equation}\label{eqn:log}
|\log(1+x)-x| \leq x^2
\end{equation}
holds for all $x \in [-1/2,1/2]$. Therefore, writing in \eqref{9d-3}, for any $M=1,\dots,d$,
\[ \log\left(\frac{\widehat{V}_{N, M}}{ \sqrt{ \vert \mathcal{L}_{N, \gamma} \vert }}+1 \right) = \log\left(\frac{\widehat{V}_{N, M}}{ \sqrt{ \vert \mathcal{L}_{N, \gamma} \vert }}+1 \right)-\frac{\widehat{V}_{N, M}}{ \sqrt{ \vert \mathcal{L}_{N, \gamma} \vert }}+\frac{\widehat{V}_{N, M}}{ \sqrt{ \vert \mathcal{L}_{N, \gamma} \vert }}\]
and combining \eqref{eqn:log} with Proposition \eqref{pp6}, we deduce
\[ 	 \sqrt{ \vert \mathcal{L}_{N, \gamma} \vert }( H-\widehat{H}_{N}) \sim - \frac{1}{2}\frac{\left( \sum_{M=1}^{d} \widehat{V}_{N, M}\log M  \right)-\left(\sum_{M=1}^{d}\widehat{V}_{N, M}\right)\left(\sum_{M=1}^{d} \log M  \right) }{ \left(\sum_{M=1}^{d} (\log M )^2 \right)-\left(\sum_{M=1}^{d} \log M  \right)^2} \]
as $N \to + \infty$.
The limit theorem (\ref{9d-2}) easily follows, thanks to Corollary \ref{cor1}. \qed

\section{Appendix}
	The basic tools from the analysis on Wiener space are presented in this section. We will focus on some elementary facts about multiple stochastic integrals. We refer to \cite{N} or \cite{NP-book} for a complete review on the topic. 

Consider ${\mathcal{H}}$ a real separable infinite-dimensional Hilbert space
with its associated inner product ${\langle
	\cdot,\cdot\rangle}_{\mathcal{H}}$, and $(B (\varphi),
\varphi\in{\mathcal{H}})$ an isonormal Gaussian process on a
probability space $(\Omega, {\mathcal{F}}, \mathbf{P})$, which is a
centered Gaussian family of random variables such that
$\mathbf{E}\left( B(\varphi) B(\psi) \right) = {\langle\varphi,
	\psi\rangle}_{{\mathcal{H}}}$ for every
$\varphi,\psi\in{\mathcal{H}}$. Denote by $I_{q} (q\geq 1)$ the $q$th multiple
stochastic integral with respect to $B$, which is an
isometry between the Hilbert space ${\mathcal{H}}^{\odot q}$
(symmetric tensor product) equipped with the scaled norm
$\sqrt{q!}\,\Vert\cdot\Vert_{{\mathcal{H}}^{\otimes q}}$ and
the Wiener chaos of order $q$, which is defined as the closed linear
span of the random variables $H_{q}(B(\varphi))$ where
$\varphi\in{\mathcal{H}},\;\Vert\varphi\Vert_{{\mathcal{H}}}=1$ and
$H_{q}$ is the Hermite polynomial of degree $q\geq 1$ defined
by :\begin{equation}\label{Hermite-poly}
	H_{q}(x)=\frac{(-1)^{q}}{q!} \exp \left( \frac{x^{2}}{2} \right) \frac{{\mathrm{d}}^{q}%
	}{{\mathrm{d}x}^{q}}\left( \exp \left(
	-\frac{x^{2}}{2}\right)\right),\;x\in \mathbb{R}.
\end{equation}
For $q=0$, 
\begin{equation}\label{27a-1}
	\mathcal{H}_{0}=\mathbb{R} \mbox{ and }I_{0}(x)=x \mbox{ for every }x\in \mathbb{R}.
\end{equation}
The isometry property of multiple integrals can be written as follows : for $p,\;q\geq
0$,\;$f\in{{\mathcal{H}}^{\otimes p}}$ and
$g\in{{\mathcal{H}}^{\otimes q}}$
\begin{equation} \mathbf{E}\Big[I_{p}(f) I_{q}(g) \Big]= \left\{
	\begin{array}{rcl}\label{iso}
		q! \langle \tilde{f},\tilde{g}
		\rangle _{{\mathcal{H}}^{\otimes q}}&&\mbox{if}\;p=q,\\
		\noalign{\vskip 2mm} 0 \quad\quad&&\mbox{otherwise,}
	\end{array}\right.
\end{equation}
where $\tilde{f}$ stands for the symmetrization of $f$. When $\mathcal{H}= L^{2}(T)$, with $T$ being an interval of $\mathbb{R}$, we have the following product formula: for $p,\;q\geq
0$,\;  $f\in{{\mathcal{H}}^{\odot p}}$ and
$g\in{{\mathcal{H}}^{\odot q}}$,  

\begin{eqnarray}\label{prod}
	I_{p}(f) I_{q}(g)&=& \sum_{r=0}^{p \wedge q} r! \binom{q}{r} \binom{p}{r}I_{p+q-2r}\left(f\otimes_{r}g\right),
\end{eqnarray}
where, for $r=0, ..., p\wedge q$, the contraction $f\otimes _{r} g$ is the function in $L ^{2}( T ^{p+q-2r}) $ given by 

\begin{equation}\label{contra}
	(f\otimes _{r} g) (t_{1},..., t_{p+q-2r})= \int_{T^{r}} f(u_{1},..., u_{r}, t_{1},..., t_{p-r}) g(u_{1},..., u_{r}, t_{p-r+1},...,t_{p+q-2r}) du_{1}...du_{r}.
\end{equation}

An useful  property of  finite sums of multiple stochastic integrals is the hypercontractivity. Namely,  for every fixed real number $p\geq 2$, there exists a universal deterministic finite constant $C_p$, such that, for any random variable $F$ of the form $F= \sum_{k=0} ^{n} I_{k}(f_{k}) $ with $f_{k}\in \mathcal{H} ^{\otimes k}$, the following inequality holds:
\begin{equation}
	\label{hyper}
	\mathbf{E}\big[\vert F \vert ^{p}\big] \leq C_{p} \left( \mathbf{E}\big[F ^{2}\big] \right) ^{\frac{p}{2}}.
\end{equation}

We denote by $D$ the Malliavin derivative operator that acts on
cylindrical random variables of the form $F=g(B(\varphi
_{1}),\ldots,B(\varphi_{n}))$, where $n\geq 1$,
$g:\mathbb{R}^n\rightarrow\mathbb{R}$ is a smooth function with
compact support and $\varphi_{i} \in {{\mathcal{H}}}$, in the following way:
\begin{equation*}
	DF=\sum_{i=1}^{n}\frac{\partial g}{\partial x_{i}}(B(\varphi _{1}),
	\ldots , B(\varphi_{n}))\varphi_{i}.
\end{equation*}
The operator $D$ is closable and it can be extended to $\mathbb{D} ^{1, 2}$ which denotes the closure of the set of cylindrical random variables with respect to the norm $\Vert \cdot\Vert _{1,2}$ defined as
\begin{equation*}
	\Vert F\Vert _{1,2} ^{2}:= \mathbf{E}\big[\vert F\vert ^{2}\big]+ \mathbf{E} \big[\Vert  DF\Vert _{\mathcal{H}} ^{2}\big]. 
\end{equation*} 
If $F=I_{p}(f)$, where $f\in \mathcal{H} ^{\odot p}$ with $\mathcal{H}= L^{2}(T)$ and $p\geq 1$, then 
$$D_{\ast}F=pI_{p-1} \left( f(\cdot, \ast)\right),$$
where $"\cdot "$ stands for $p-1$ variables.

The pseudo inverse $ (-L) ^{-1}$ of the Ornstein-Uhlenbeck operator $L$ is defined, for $F=I_{p}(f)$ with $f\in \mathcal{H} ^{\odot p}$ and $p\geq 1$, by 
\begin{equation*}
	(-L) ^{-1}F= \frac{1}{p} I_{p} (f).
\end{equation*}

At last notice that in our work, we have $\mathcal{H}= L^2(\mathbb{R})$ while the role of the isonormal process $(B(\varphi), \varphi \in \mathcal{H}) $ is played by  the usual Wiener integral on $L^2(\mathbb{R})$ associated with the Brownian motion $(B(y), y \in \mathbb{R})$. In this case, we can provide an explicit formula for the multiple integral $I_q$, ($q\geq 1$). Indeed, if $f$ is a symmetric function of the form
\begin{equation}\label{eqn:Wiener-Itôintegral}
f = \sum_{j_1,\ldots,j_q=1}^n a_{j_1,\ldots,j_q} \mathbbm{1}_{[s_{j_1},t_{j_1})} \otimes \cdots \otimes \mathbbm{1}_{[s_{j_q},t_{j_q})},
\end{equation}
where, $\otimes$ stands for the tensor product, $a_{j_1,\ldots,j_d}$ are such that, for all permutation $\sigma$, $a_{\sigma(j_1),\ldots,\sigma_(j_q)}=a_{j_1,\ldots,j_q}$ and $a_{j_1,\ldots,j_q}=0$ as soon as two indices $j_1,\ldots,j_q$ are equal and, for all $1 \leq \ell \neq \ell' \leq q$, $[s_{j_\ell},t_{j_\ell}) \cap [s_{j_{\ell'}},t_{j_{\ell'}}) = \emptyset$, then
\begin{equation}\label{eqn:Wiener-Itôintegral2}
I_q(f) := \sum_{j_1,\ldots,j_q=1}^n a_{j_1,\ldots,j_q} (B(t_{j_1})-B( s_{j_1})) \times \ldots (B(t_{j_q})-B( s_{j_q})).
\end{equation} 
It is straightforward that this last random variable belongs to $L^2(\Omega)$. For a general symmetric $f \in L^2(\mathbb{R}^q)$, $I_q(f)$ is then defined using the density of functions of the form \eqref{eqn:Wiener-Itôintegral} within the set of symmetric square integrable function and by checking that the corresponding random variables \eqref{eqn:Wiener-Itôintegral2} converge in $L^2(\Omega)$.

Our main tool to prove the asymptotic normality  for random vectors is the following theorem (Theorem 6.1.1 in \cite{NP-book}). 

\begin{theorem}
	\label{clt}
	Let $m\geq 1$ be an integer number and consider a $m$-dimensional random vector $ F=( F_{1},...,F_{m})$. Assume $ F_{i} \in \mathbb{D} ^{1,4}$ for every $i=1,...m$. Let $C \in M _{m}(\mathbb{R})$ be a symmetric and positive definite matrix and let $Z\sim N(0, C)$. Then
	\begin{equation}
		d_{W}(F, Z) \leq C \sqrt{ \sum_{i,j=1} ^{m} \mathbf{E} \Big[(C_{i,j}-\langle DF_{i}, D(-L) ^{-1} F_{j} \rangle ) ^{2}\Big] }.
	\end{equation}

\end{theorem}

Let us finally mention the fact that, if $F_i \in \mathbb{D}^{1,2}$ for every $i=1,...m$, as a consequence of \cite[Theorem 2.9.1]{NP-book},

\begin{equation}\label{cltbis}
\mathbf{E}[\langle DF_{i}, D(-L) ^{-1} F_{j} \rangle ] = \mathbf{E}[F_iF_j].
\end{equation}

\bibliography{biblio}{}
\bibliographystyle{plain}

\end{document}